\documentclass[10pt,leqno]{amsart}

\usepackage{amsmath,amsthm,amscd,amssymb,amsfonts, amsbsy}
\usepackage{latexsym}
\usepackage{txfonts}
\usepackage{exscale}
\usepackage{mathrsfs}
\usepackage{esint} % added to replace the $\nint$ with $\fint$
\usepackage{enumitem} 
\usepackage[pagewise]{lineno}
\usepackage{color}
\usepackage[useregional]{datetime2}
\usepackage[
  hmarginratio={1:1},     % equal left and right margins
  vmarginratio={1:1},     % equal top and bottom margins
  textwidth=400pt,			  % new text width
	textheight=580pt,        % new text height
  heightrounded,          % always useful
]{geometry}

\usepackage[utf8]{inputenc}
\usepackage[colorlinks,citecolor=blue,hypertexnames=false]{hyperref}

\usepackage{appendix}

\allowdisplaybreaks

\numberwithin{equation}{section}

\theoremstyle{plain}
\newtheorem{theorem}[equation]{Theorem}
\newtheorem{prop}[equation]{Proposition}
\newtheorem{corollary}[equation]{Corollary}
\newtheorem{lemma}[equation]{Lemma}

\theoremstyle{definition}
\newtheorem{defn}[equation]{Definition}

\theoremstyle{remark}

\newtheorem{remark}[equation]{Remark}

\numberwithin{equation}{section}

\newcommand{\RR}{{\mathbb{R}}}
\newcommand{\NN}{{\mathbb{N}}}

\newcommand{\rn}{\mathbb{R}^n}

\def\XXint#1#2#3{{\setbox0=\hbox{$#1{#2#3}{\int}$}
     \vcenter{\hbox{$#2#3$}}\kern-.5\wd0}}

\reversemarginpar

\begin{document}
\allowdisplaybreaks

\title[]{On the generalized parabolic Hardy-H\'enon  equation: Existence, blow-up, self-similarity\\ and large-time asymptotic behaviour}

\author[]{Gael Diebou Yomgne}

\address{Gael Diebou Yomgne,
Mathematisches Institut, Rheinische Friedrich-Wilhelms-Universit\"{a}t Bonn\\
Endenicher Allee 60,  53115, Germany} \email{gaeldieb@math.uni-bonn.de}

\thanks{The author is supported by the DAAD through the program \textquotedblleft Graduate School Scholarship Programme, 2018   number 57395813\textquotedblright and by  Hausdorff Center for Mathematics at Bonn}
\date{\today}

\begin{abstract}
This paper deals with the Cauchy problem for the Hardy-H\'{e}non equation (and its fractional analogue). Local well-posedness for initial data in the class of continuous functions with slow decay at infinity is investigated.
Small data (in critical weak-Lebesgue space) global well-posedness is obtained in $C_{b}([0,\infty);L^{q_c,\infty}(\RR^{n}))$. As a direct consequence, global existence for data in strong critical Lebesgue $L^{q_c}(\RR^{n})$ follows under a smallness condition while uniqueness is unconditional. Besides, we prove the existence of self-similar solutions and examine the long time behavior of globally defined solutions. The zero solution $u\equiv 0$ is shown to be asymptotically stable in $L^{q_c}(\RR^{n})$ -- it is the only self-similar solution which is initially small in $L^{q_c}(\RR^{n})$. Moreover, blow-up results are obtained under mild assumptions on the initial data and the corresponding Fujita critical exponent is found.

\end{abstract}

\maketitle

%\tableofcontents[hideallsubsections]
\section{Introduction}
While a great deal of works on the Cauchy problem for the reaction diffusion equation $\partial_tu-\Delta u=|u|^{p-1}u$ in various geometries \cite{Caze-Weiss,Lee-Ni,Giga,Weissler1,Snou-Taya-Weiss} and the hyper-dissipative equation (see e.g. \cite{Gala-Poza,Ferreira-Villa,Caristi-Mitie,Grunau} and references within) already exist in the literature, several aspects of solutions to their analogues with singular nonlinearity are yet to be investigated.
In this paper, we are interested in the Cauchy problem for the Hardy-H\'{e}non equation
\begin{equation}\label{eq:main-eq}
u_t+(-\Delta)^{m}u=|x|^{-\alpha}F(u) \ \quad\mbox{for }
(x,t)\in\mathbb{R}^{n}\times(0,\infty);\quad u(0)=u_0 \hspace{0.2cm}\mbox{on } x\in\mathbb{R}^{n}
\end{equation}
where  $m\in (0,1)\cup \mathbb{N}$, $\alpha\in \RR$  and $(-\Delta)^{m}$ is the power $m >0$ of the minus Laplace operator $(-\Delta)$ which transforms any Schwartz function $\phi\in \mathcal{S}(\rn)$ according to $\widehat{(-\Delta)^{m}\phi}(\xi)=|\xi|^{2m}\widehat{\phi}(\xi)$ where $\widehat{\phi}$ denotes the Fourier transform of $\phi$. The nonlinearity $F$ is such that
\begin{equation}\label{eq:cond-on-F}
F(0)=0,\quad|F(a)-F(b)|\leq C|a-b|(|a|^{p-1}+|b|^{p-1}) 
\end{equation}
for some constant $C>0$ and for all $a,b\in \RR$, $p>1$.

 Probably the well-posedness of the Cauchy problem for the Hardy-H\'{e}non equation (i.e. Eq. \eqref{eq:main-eq} for $m=1$) was first studied by Wang \cite{Wang} and has recently attracted a lot of attention. We shall mention the work of the authors in \cite{Slim-Tay-Weiss} whereby local and global existence questions have been addressed (for $F(u)=u|u|^{p-1}$ and $\alpha>0$). Indeed, for initial data in $C_{0}(\RR^{n})$ (the set of continuous functions vanishing at infinity) and in $L^{q}(\RR^{n})$, $q>\frac{n(p-1)}{2-\alpha}$ local mild solutions were constructed -- their blow up rate was also determined. As for global-in-time solutions, it was proved that they arise from small initial data in $L^{\frac{n(p-1)}{2-\alpha}}(\RR^{n})$ and are small in $C([0,\infty);L^{\frac{n(p-1)}{2-\alpha}}(\RR^{n}))$. On the other hand if the prescribed data decay faster at infinity than the singularity $\omega(x)|x|^{-\frac{2-\alpha}{p-1}}$, $\omega\in L^{\infty}(\RR^{n})$ with small norm, then the corresponding problem has a global mild solution which at large times is asymptotic to a self-similar solution in the $L^{\infty}$-topology. They further investigated the asymptotics of solutions to a modified model involving a potential multiplying the singularity in problem \ref{eq:main-eq}. In the same vein, the well-posedness and asymptotic stability of Eq. \eqref{eq:main-eq} in Besov space was considered in \cite{Noboru}. Other interesting results in the local theory have recently been obtained in \cite{Tayachi}.

Here, we carry out the analysis of Problem \ref{eq:main-eq} in a different functional setting.  If $u(x,t)$ solves Eq. \eqref{eq:main-eq} in a classical sense, then for each $\sigma>0$, $u_{\sigma}(x,t)=\sigma^{\frac{2m-\alpha}{p-1}}u(\sigma x,\sigma^{2}t)$ is another solution with initial data $u_{\sigma,0}(x)=\sigma^{\frac{2m-\alpha}{p-1}}u_0(\sigma x)$ provided, of course that $F$ satisfies the homogeneity property $F(u(\sigma x,\sigma^{2m}t))=\sigma^{-\frac{(2m-\alpha)p}{p-1}}F(u_{\sigma})(x,t)$.  Observe that as initial data, $f(x)=\varepsilon_0|x|^{\frac{-(2m-\alpha)}{p-1}}$ has such a property -- we anticipate that it generates a global solution, that is, a self-similar solution whenever $\varepsilon_0$ is sufficiently small. However, $f\notin L^{q}(\RR^{n})$ for every $q>0$ while $f\in L^{q,r}(\RR^{n})$ for the unique choice $q=q_c=\frac{n(p-1)}{2m-\alpha}$ and $r=\infty$. This suggests the Marcinkiewicz space $L^{q_c,\infty}(\RR^{n})$ as a natural setting for the analysis of the Cauchy problem for Eq. \eqref{eq:main-eq}, notably the study of self-similar solutions. This practice goes back (at least) to the work of Giga \& Miyakawa \cite{Giga-Miya} in the context of Navier-Stokes equations and has inspired several works in the literature. For small data in $L^{q_c,\infty}(\RR^{n})$, we establish (see Theorem \ref{thm:GWP-Lp}) existence of a small global-in-time solution in the space  $C_b([0,\infty);L^{q_c,\infty}(\mathbb{R}^n))$. In Theorem \ref{thm:uniqueness}, it is shown that this smallness condition can be dropped so that unconditional uniqueness holds whenever $u_0$ belongs to the critical strong Lebesgue space, hence improving \cite[Theorem 1.3 (i)]{Slim-Tay-Weiss} and other earlier uniqueness results. In the process, we essentially rely on smoothing effect type estimates for the singular polyharmonic heat semigroup in weak-Lebesgue spaces (see Proposition \ref{prop:smoothing-effect}), as well as analogues of Yamazaki's bounds \cite{Yama} (cf. Lemma \ref{lem:nonlinearbound}) pertaining to some nonlinear operator.
Self-similar solutions are also investigated and are relevant in analyzing the long time behavior of globally defined solutions. Actually, they appear to be attractors of global solutions whose initial values are suitable perturbations of homogeneous functions lying in weak-$L^{q_c}(\RR^{n})$ thus giving a characterization of the basin of attraction of each attractor. In the strong critical Lebesgue framework, the global solution decays to zero at large times. On the other hand, initial data decaying faster than $|x|^{-\gamma}$, $\gamma=\frac{2m-\alpha}{p-1}$ produce global-in-time solutions, see Theorem \ref{thm:GWP-decay}. Interestingly, for positive data and nonlinearity, this decay assumption seems sharp in a certain sense as pointed out in Remark \ref{rmk:GWP}. In addition, if the nonlinearity is positive, then provided $1<p<p_F=1+\frac{2m-\alpha}{n}$ and $u_0\in L^1(\RR^n)$ has a positive integral, Eq. \eqref{eq:main-eq} has no solution  so that the critical Fujita exponent for Eq. \eqref{eq:main-eq} is $p_F$, see Theorem  \ref{thm:nonexistence}. The proof of this blow-up result employs the test function method \cite{Miti-Poho} which requires that we consider solutions in the sense of Definition \ref{defn:Weak-sol} below. Thus, the first step towards the proof consists of showing that mild solutions are indeed weak solutions under reasonable assumptions on the initial data.  To be consistent with the terminology used above, we record the following
\begin{defn}\label{defn:Mild-sol}Let $\mathcal{X}$ be a Banach space of real-valued functions. We say that a function $u$ in $C\big([0,T);\mathcal{X}\big)$, $T>0$ is a global-in-time (resp. local-in-time) mild solution of the Cauchy problem for Eq. \eqref{eq:main-eq} if $u$ satisfies the integral equation
\begin{equation}\label{eq:mild-solution}
u(t)=e^{-t(-\Delta)^m}u_0+\int_{0}^{t}e^{-(t-\tau)(-\Delta)^m}|\cdot|^{-\alpha}F(u(\tau))d\tau
\end{equation}	
for every $t>0$ (i.e. $T=\infty$) (resp. for $t\in (0,T)$ with $T:=T(u_0)<\infty$) and $u(t)$ converges to $u_0$ in a suitable sense as $t\rightarrow 0^+.$
\end{defn}	
 $C\big([0,T);\mathcal{X}\big)$ denotes the space of continuous functions of $t\in [0,T)$ with values in $\mathcal{X}$. Often, we will use $C_b\big([0,T);\mathcal{X}\big)$ to denote the space of bounded and continuous functions from $[0,T)$ onto $\mathcal{X}$ endowed with the space-time norm 
 \[[u]_{\mathcal{X}}:=\sup_{t>0}\|u(t)\|_{\mathcal{X}}.\]
For those spaces frequently used in our analysis, namely Lebesgue spaces,  $C_{\Lambda}(\RR^n)$ (see below for its definition) the convergence in the above definition is taken with respect to the strong topology in the corresponding space. In case $\mathcal{X}$ is a dual of a Banach space, e.g. the weak-Lebesgue space $L^{q,\infty}(\RR^n)$, convergence as $t\rightarrow 0^{+}$ should be understood in the weak-$\star$ sense.\\
Another notion of solution we consider in the sequel is the so-called weak solution, by which we mean the following: 
\begin{defn}\label{defn:Weak-sol}
Let $F(u)$ as in \eqref{eq:cond-on-F} and take $\alpha\in (0,\min\{2m,n\})$, $m\in \mathbb{N}$. A function $u: \RR^{n}\times(0,\infty)\rightarrow \RR$ is said to be a global weak solution of Eq. \eqref{eq:main-eq} if it obeys 
\begin{enumerate}
\item $u_0\in L^{1}_{loc}(\RR^{n})$,\quad
	 $|x|^{-\alpha}F(u)\in L^{1}_{loc}( \RR^{n}\times(0,\infty))$	
\item $\displaystyle\int_{\RR^{n}\times(0,\infty)}u(-\partial_{t}\psi+(-\Delta)^{m}\psi)dxdt=\int_{\RR^{n}}u_0(x)\psi(x,0)dx+\int_{\RR^{n}\times(0,\infty)}|x|^{-\alpha}|u|^{p}\psi dxdt$
\end{enumerate}
for all $\psi\in C^{\infty}_{0}( \RR^{n}\times\mathbb{R})$.
\end{defn}

Throughout this paper, unless otherwise specified, the dimension will always be taken larger or equal to $1$. The generic constants $C$ and $c$ which will be appearing often in chain of estimates may differ from one line to another and their dependence on other parameters will be highlighted whenever necessary.

\section{Main results}\label{S:main-results}
Our first result pertains to the local well-posedness subject to initial data decaying at infinity.

\begin{defn}
Consider the function $\Lambda(x)=(1+|x|)^{\frac{-\alpha}{p-1}}$, ($p>1$, $\alpha<0$). We say that a function $\phi$ belongs to $C_{\Lambda}(\RR^{n})$ if $\phi$ is continuous on $\RR^{n}$ and satisfies \[\|\phi\|_{C_{\Lambda}}:=\|\phi\|_{\Lambda}=\|\Lambda \phi\|_{L^{\infty}(\RR^n)}<\infty.\]
\end{defn}
\begin{theorem}\label{thm:LWP-positive-alpha}
Assume that $\alpha<0$ and $n\geq1$. For any $u_0\in C_{\Lambda}(\RR^{n})$, the Cauchy problem \ref{eq:main-eq} has a unique classical solution $u$ on $\RR^{n}\times [0,T_{\infty})$ with either $T_{\infty}=\infty$ or $T_{\infty}<\infty$ and $\displaystyle \lim_{t\rightarrow T_{\infty}}\|u\|_{\Lambda}=\infty$.	
\end{theorem} 

The existence and uniqueness statement for globally defined solutions of Eq. \eqref{eq:main-eq} is given below.
\begin{theorem}\label{thm:GWP-Lp}
Let $0<\alpha<2m<n$. Assume that F satisfies \eqref{eq:cond-on-F}. We have the following assertions:
\begin{enumerate}[label={($\textbf{B}_{\arabic*})$}]
\item \label{GWP-P1} There exist two constants $\varepsilon>0$ and $\vartheta:=\vartheta(\varepsilon)>0$ such that for any $u_0\in L^{q_c,\infty}(\RR^{n})$ with $\|u_0\|_{L^{q_c,\infty}}<\varepsilon$, there exits a global-in-time solution $u$ of \eqref{eq:mild-solution} in $C_{b}([0,\infty);L^{q_c,\infty}(\RR^{n}))$ which is unique in the closed ball \[\overline{B}_{2\vartheta}=\{u\in C_{b}([0,\infty);L^{q_c,\infty}(\RR^{n})): [u]_{L^{q_c,\infty}}\leq 2\vartheta\},\] provided $\frac{n-\alpha}{n-2m}<p<\infty$.
\item \label{GWP-P2}  Given $u_0\in (L^{r,\infty}\cap L^{q_c,\infty})(\RR^{n})$ with $r'\in (1,\frac{n}{2m-\alpha})$, there exists $\tilde{\varepsilon}:=\tilde{\varepsilon}(r)>0$, $\tilde{\varepsilon}<\varepsilon$ such that if $\|u_0\|_{L^{q_c,\infty}}<\tilde{\varepsilon}$, then a version of the above conclusion is true in the space $C_{b}([0,\infty);L^{r,\infty}(\RR^{n}))$. 
\item\label{GWP-P3} Let $0<2m+\alpha<n$ and $\frac{n-\alpha}{n-2m}<p<\frac{2m}{\alpha}$. Assume further that $q_c<q<\frac{np(p-1)}{2m-\alpha p}$. Let $\varepsilon>0$ as in \ref{GWP-P1} and put $\delta=\frac{n}{2m}(\frac{1}{q_c}-\frac{1}{q})$. There exists $\varepsilon_q\in (0,\varepsilon)$ such that for all $u_0\in  L^{q_c,\infty}(\RR^{n})$ with $\|u_0\|_{L^{q_c,\infty}}<\varepsilon_q$, there is exactly one solution $u$ of Eq. \eqref{eq:mild-solution} lying in the space 
$$S^{\delta}_{q}=\{u\in C_{b}([0,\infty);L^{q_c,\infty}(\RR^{n})):t^{\delta}u\in C_{b}((0,\infty);L^{q,\infty}(\RR^{n}))\}.$$ 
\end{enumerate}
\end{theorem}
Our result via \ref{GWP-P1} shows in particular that global-in-time solutions which are initially (at time $t=0$) singular exist. However, we do not know whether those solutions enjoy a higher regularity. In case $m=1$, the answer is known to be negative for bounded continuous initial data, see \cite{Wang}. The solution constructed in the above theorem is unique only in a small ball of the underlying solution space. If the data is prescribed in the critical strong Lebesgue space, then uniqueness holds for free.
 \begin{theorem}[\textbf{Unconditional uniqueness}]\label{thm:uniqueness}
Assume that $p\in (\frac{n-\alpha}{n-2m},\infty)$ with $\alpha\in (0,2m)$ with $2m<n$. Consider $u$, $v$ two global solutions of Eq. \eqref{eq:mild-solution} in   $C_{b}([0,\infty);L^{q_c}(\RR^{n}))$. Then $u=v$ a.e. in $\RR^{n}\times (0,\infty)$. 
\end{theorem}
This theorem is stronger than the uniqueness results obtained in \cite[Corollary 1]{Slimene} and \cite{Slim-Tay-Weiss} for the special case $m=1$ via the two-norms approach and under a smallness condition on the data. A similar result was obtained in \cite{Ferreira-Villa} in the non-singular case $\alpha=0$. A scaling analysis predicts that well-posedness should hold whenever $p\in (\frac{2m-\alpha}{n}+1,\infty)$. However, Theorem \ref{thm:GWP-Lp} is valid for $p$ in a sub-range of the latter interval due the choice of the functional framework. This gap is close in our next result.

\begin{theorem}\label{thm:GWP-decay}	
Let $0<\alpha<\min{\{2m,n\}}$ and assume $q_c>1$. Let $u_0\in L^{1}_{loc}(\RR^{n})$ with the property that $|u_0(x)|\leq c_0|x|^{-\frac{2m-\alpha}{p-1}}$ for all $x\in \RR^{n}\setminus\{0\}$ and $c_0>0$ small enough. Then,  Eq. \eqref{eq:mild-solution} is  globally well-posed.
\end{theorem}
This improves \cite[Theorem 2.1]{Caristi-Mitie} since no sign condition is imposed on $u_0$. A version of the above result for $m=1$ can be found in \cite{Slim-Tay-Weiss}. A natural question to ask is whether this fast decay assumption on the data is optimal. It turns out that the answer is affirmative, at least for positive nonlinearity.
\begin{theorem}[\textbf{Nonexistence}]\label{thm:nonexistence}
Let $m\in \mathbb{N}$, $\alpha\in (0,\min\{n,2m\})$ and assume that $1<p<p_F$. Put $F(u)=|u|^{p}$. Given $u_0\in L^{1}(\RR^{n})\cap C_0(\rn)$ with $\displaystyle\int_{\RR^{n}}u_0(x)dx>0$, there exists no mild solution of Eq. \eqref{eq:main-eq}. 
\end{theorem}
The proof of this result relies on the fact that  mild solutions  are weak solutions of \eqref{eq:main-eq} in the sense of Definition \ref{defn:Weak-sol}.
\begin{remark}\label{rmk:GWP}A few observations are in order:
\begin{enumerate}[label=(\roman*)] 
\item Theorem \ref{thm:GWP-Lp} remains valid for $m\in (0,1)$ i.e. the fractional version of the Hardy-H\'{e}non equation. The integral kernel of the semigroup generated by the nonlocal operator $-(-\Delta)^{m}$, $0<m<1$ has an algebraic decay at infinity (see e.g. \cite{Jacob}), so that Proposition \ref{prop:smoothing-effect} below, heavily used in our analysis naturally extends to this case. 
\item As revealed by the proof of Theorem \ref{thm:nonexistence}, one can show in parallel that if $u_0$ is nonnegative a.e. on $\RR^{n}$ and there exists  $\beta\in [0,\frac{2m-\alpha}{p-1})$ such that $$0<K\leq \displaystyle\liminf_{|x|\rightarrow\infty}(1+|x|^{\beta})u_0(x),$$ then Problem \ref{eq:main-eq} has no weak solution (see \cite{Caristi-Mitie} for a similar result in the nonsingular case). Combining Theorems \ref{thm:GWP-decay} and \ref{thm:nonexistence}, we find that $p_F=1+\dfrac{2m-\alpha}{n}$ is the Fujita critical exponent for Eq. \eqref{eq:main-eq}.
\end{enumerate}	
	   
\end{remark}

\section{Preliminaries and auxiliary results}\label{s3}
We establish in this section relevant estimates on the singular polyharmonic heat semigroup (and its kernel) in the framework of Lorentz spaces which play a pivotal role in the proof of our main results.  
\subsection{Lorentz spaces and their properties}
We briefly recall those properties of Lorentz spaces according to their usefulness later in the paper. For a systematic description, we do refer the reader to \cite{Bennett-Sharpley}. 
Denote by $|E|$ the Lebesgue measure of a measurable set $E\subset \rn$. Call $f^{\ast}:[0,\infty)\rightarrow [0,\infty)$ the decreasing rearrangement of a measurable function $f:\rn\rightarrow \RR$ which is defined as $f^{\ast}(\lambda)=\inf\{\tau>0:|\{|f|> \tau\}|\leq \lambda\}$ and let $f^{\ast\ast}:(0,\infty)\rightarrow [0,\infty)$ be the associated maximal function given by $f^{\ast\ast}(\tau)=\frac{1}{\tau}\int_{0}^{\tau}f^{\ast}(s)ds$.  
\begin{defn}\label{defn:Lorentz-spaces}
Let $1<p<\infty$, $1\leq q\leq \infty$, the Lorentz space $L^{p,q}(\rn)$ is defined as
\[L^{p,q}(\rn)=\bigg\{f:\rn\rightarrow \RR \hspace{0.1cm}\text{measurable} : \|f\|^{q}_{L^{p,q}}=\int_{0}^{\infty}[s^{1/p}f^{\ast\ast}(s)]^{q}\frac{ds}{s}<\infty\bigg\},
\]
and
\[L^{p,\infty}(\rn)=\bigg\{f:\rn\rightarrow \RR \hspace{0.1cm}\text{measurable} : \|f\|_{L^{p,\infty}}=\sup_{0<s<\infty}s^{1/p}f^{\ast\ast}(s)<\infty\bigg\}.
\]
\end{defn}
These spaces increase (with continuous embedding) with respect to the second exponent, that is, $L^{p,q}(\rn)\subseteq L^{p,r}(\rn)$ for any $q\leq r\leq \infty$ and share a few common properties with Lebesgue spaces $L^{p}(\rn)=L^{p,p}(\rn)$ such as the scaling law, $\|f(\sigma\cdot)\|_{L^{p,q}}=\sigma^{-\frac{n}{p}}\|f\|_{L^{p,q}}$ for all $\sigma>0$ as well as the following identity $\||f|^{l}\|_{L^{p,q}}=\|f\|^{l}_{L^{lp,lq}}$ for all $l>0$. Also the K\"{o}the dual space (or associate space) of $L^{p,q}(\rn)$, $1<p,q<\infty$ is $L^{p',q'}(\rn)$ where $p',q'$ are conjugate exponents of $p$ and $q$ respectively while for $p\in (1,\infty)$ and $q\in (0,1]$, we have $(L^{p,q}(\rn))'=L^{p',\infty}(\rn)$. Note that Definition \ref{defn:Lorentz-spaces} does not cover the case where $0<p,q<1$. In fact for this choice, $L^{p,q}(\RR^{n})$ is a quasi-Banach space defined exactly as in Definition \ref{defn:Lorentz-spaces} with $f^{\ast\ast}$ replaced by $f^{\ast}$.
 A version of H\"{o}lder's inequality and Young's inequality hold in the setting of Lorentz spaces (see \cite{Oneil}). 
\begin{lemma}[Generalized H\"{o}lder's inequality]\label{lem:Holder ineq}
 Let $1<p_1,p_2<\infty$ and $1\leq q_1,q,q_2\leq \infty$ such that $\frac{1}{q}\leq  \frac{1}{q_1}+\frac{1}{q_2}$. For all $f\in L^{p_1,q_1}(\rn)$ and $g\in L^{p_2,q_2}(\rn)$, we have $fg\in  L^{p,q}(\rn)$ and it holds that 
 \[\|fg\|_{ L^{p,q}}\leq C\|f\|_{L^{p_1,q_1}}\|g\|_{L^{p_2,q_2}}\]
provided $\frac{1}{p}= \frac{1}{p_1}+\frac{1}{p_2}$. 
\end{lemma}
\begin{lemma}[Generalized Young's inequality]\label{lem:Young ineq}
Let $1<p_1,p_2<\infty$ and $1\leq q_1,q_2\leq \infty$.
Assume $\frac{1}{p_1}+ \frac{1}{p_2}>1$, $\frac{1}{p}+1=\frac{1}{p_1}+\frac{1}{p_2}$, $p\in (1,\infty]$ and let $f\in L^{p_1,q_1}(\rn)$, $g\in L^{p_2,q_2}(\rn)$. We have that $f\ast g\in L^{p,q}(\rn)$ as long as $\frac{1}{q}\leq \frac{1}{q_1}+\frac{1}{q_2}$ and there exists a constant $C>0$ with
\[\|f\ast g\|_{ L^{p,q}}\leq C\|f\|_{L^{p_1,q_1}}\|g\|_{L^{p_2,q_2}}.\] 
Moreover, in case $\frac{1}{p_1}+ \frac{1}{p_2}=1$ and $1\leq \frac{1}{q_1}+ \frac{1}{q_2}$, then the above conclusion holds for $p=q=\infty$, that is, $f\ast g\in L^{\infty}(\rn)$ and there exists a constant $c>0$ depending on the exponents such that
\[\|f\ast g\|_{L^{\infty}}\leq c\|f\|_{L^{p_1,q_1}}\|g\|_{L^{p_2,q_2}}.
\]
\end{lemma}  
\subsection{Polyharmonic linear heat equation and key estimates}
Consider the linear parabolic equation
\[\partial_{t}u+(-\Delta)^{m}u=0\hspace{0.1cm}\text{in}\hspace{0.1cm}\RR^{n}\times (0,\infty)
\] 
subject to the initial data $u(x,0)=u_{0}(x)$, $x\in \RR^{n}$.
A standard Fourier analysis on the above equation shows that it has a smooth and radial fundamental solution $g_m$ which reads 
\begin{equation*}
g_{m}(x,t)=t^{-\frac{n}{2m}}g(t^{-\frac{1}{2m}}x,1),\hspace{0.2cm}\widehat{g_{m}}(\xi,t)=e^{-|\xi|^{2m}t}, \hspace{0.2cm}t>0. 
\end{equation*}
Hence, the solution $u$ of the Cauchy problem may be formally realized via convolution by 
$$u(x,t)=e^{-t(-\Delta)^{m}}u_0(x)=\left(g_{m}(\cdot,t)\ast u_0\right)(x)$$ whenever this representation makes sense (e.g. when $u_0\in C_{c}(\mathbb{R}^{n}))$. Observe that since $e^{-|\xi|^{2m}}$ is globally integrable, the rescaled kernel $g$ is bounded continuous on $\RR^{n}$ and vanishes as $|x|\rightarrow \infty$. In addition, we have:
\begin{lemma}\label{lem:exp-est-g}
	Let $g$ be as above. Given $\eta>0$, there exits a positive constant $C:=C(\eta,n)$ such that for $\lambda\in (0,2]$ and $y\in \mathbb{R}^{n}$, 
	\begin{equation}\label{ineq:exp-int-est}
	\int_{\mathbb{R}^{n}}|g(x)|(1+|y-\lambda x|)^{-\eta}dx\leq C(1+|y|)^{-\eta}.
	\end{equation}
\end{lemma}
\begin{proof}[Proof of Lemma \ref{lem:exp-est-g}]
	Recall that $g(x)=(2\pi)^{n/2}\displaystyle\int_{\mathbb{R}^{n}}e^{ix\cdot \xi-|\xi|^{2m}}d\xi$.
	We start with the following.\\ 
	\textbf{Claim.} Fix $N\geq 0$ integer, there exists a constant $c>0$ depending on $N$ and $n$ only such that there holds the pointwise bound  
	\begin{equation}\label{ineq:pointwise-on-g}
	|g(x)|\leq c(1+|x|)^{-N}\quad \forall \hspace{0.1cm} x\in\mathbb{R}^{n}.
	\end{equation}
	Let us momentarily defer the proof of the latter and observe that it allows us to write  for any $\eta>0$,  
	\begin{align*}
	\int_{\mathbb{R}^{n}}|g(x)|(1+|y-\lambda x|)^{-\eta}dx &\leq c\int_{\mathbb{R}^{n}}(1+|x|)^{-N}(1+|y-\lambda x|)^{-\eta}dx\\
	&\leq c(L_1+L_2)
	\end{align*} 
	where $$ L_1=\int_{\{x\in \mathbb{R}^{n}:|y-\lambda x|\leq \frac{|y|}{2}\}}(1+|x|)^{-N}(1+|y-\lambda x|)^{-\eta}dx$$
	and $$L_2=\int_{\{x\in \mathbb{R}^{n}:|y-\lambda x|>\frac{|y|}{2}\}}(1+|x|)^{-N}(1+|y-\lambda x|)^{-\eta}dx.$$
	Given that (\ref{ineq:pointwise-on-g}) is true for every $N>0$ (so in particular for $N=n+\lfloor\eta\rfloor+1$ where $\lfloor\eta\rfloor$ stands for the largest integer not exceeding $\eta$) , we have
	\begin{align*}
	L_1&\leq C\int_{\{x\in \mathbb{R}^{n}:|y-\lambda x|\leq \frac{|y|}{2}\}}(1+|x|)^{-(\lfloor\eta\rfloor+n+1)}dx\\
	&\leq C\int_{\frac{2|y|}{\lambda}}^{\infty}r^{n-1}(1+r)^{-(\lfloor\eta\rfloor+n+1)}dr\\
	&\leq C\bigg(1+\frac{2}{\lambda}|y|\bigg)^{-\eta}\\
	&\leq C(n,\eta)(1+|y|)^{-\eta}
	\end{align*}
	since $0<\lambda\leq 2$. To estimate $L_2$, we use the inequality $|y|+1\leq 2(|y-\lambda x|+1)$ and proceed as follows,
	\begin{align*}
	L_2&\leq \int_{\{x\in \mathbb{R}^{n}:|y-\lambda x|> \frac{|y|}{2}\}}(1+|x|)^{-N}(1+|y|)^{-\eta}dx\\
	&\leq 2^{\eta}(1+|y|)^{-\eta}\int_{ \mathbb{R}^{n}}(1+|x|)^{-N}dx\\
	& \leq C(1+|y|)^{-\eta}.
	\end{align*}
	We finish the proof with a justification of the Claim (\ref{ineq:pointwise-on-g}). Since $g$ is radial, it uniquely solves an ODE whose asymptotic analysis leads to a pointwise decay estimate of exponential type 
	\[ |g(x)|\leq Ce^{-d|x|^{L}}, \hspace{0.2cm}L=\frac{2m}{2m-1},\hspace{0.2cm} x\in \RR^{n}\]
for some constant $C>0$ and $d>0$ depending on $m$ and $n$, see \cite[Proposition 2.1]{Gala-Poza}. On the other hand, the function $e^{-|\xi|^{2m}}$ is smooth for $m\geq 1$ so that by a well-known property of the Fourier transform, $g(x)$ decays rapidly to zero at infinity. 
\end{proof}
\begin{remark}
Alternatively, one may use the stationary phase method, via the study of the asymptotic behavior of $g(x)$ as $|x|\rightarrow\infty$ to obtain pointwise decay bounds which may differ from the one mentioned above, see for instance appendix A in \cite{Koch-Lamm} for the special situation $m=2$. 
\end{remark}
Denoting by $E_{m,\alpha}(t)$ the singular polyharmonic heat kernel $e^{-t(-\Delta)^{m}}|\cdot|^{-\alpha}$, we have the following boundedness properties. 
\begin{prop}\label{prop:smoothing-effect}
Let $\alpha\in (0,\min\{2m,n\})$ and $1<p_1,p_2\leq \infty$ and $q\in [1,\infty]$ such that
\begin{equation}\label{ineq:condp-q}
p_1>\frac{n}{n-\alpha}, \quad p_2>\frac{np_1}{n+\alpha p_1}.
\end{equation}
Then 
\begin{enumerate}[label={($\textbf{A}_{\arabic*})$}]
		\item \label{part a}$E_{m,\alpha}(t)$ maps continuously $L^{p_1}(\mathbb{R}^{n})$ into $C_0(\mathbb{R}^{n})$, the set of continuous functions on $\mathbb{R}^{n}$ vanishing at infinity.
		\item\label{part b} For $\beta\in \NN^{n}_0$, there exists a positive constant $C:=C(n,p_1,p_2,\alpha,\beta,m)$ such that 
\begin{equation}\label{eq:Lorentz-smoothing-effect}
		\|\partial^{\beta}E_{m,\alpha}(t)f\|_{L^{p_2,q}}\leq Ct^{-\frac{n}{2m}\big(\frac{1}{p_1}-\frac{1}{p_2}\big)-\frac{\alpha}{2m}-\frac{|\beta|}{2m}}\|f\|_{L^{p_1,\infty}}\quad \text{for all}\quad t>0,\hspace{0.1cm}f\in L^{p_1,\infty}(\mathbb{R}^{n}).
		\end{equation}
		In particular,  $\partial^{\beta}E_{m,\alpha}(t): L^{p_1}(\mathbb{R}^{n})\longrightarrow L^{p_2}(\mathbb{R}^{n})$ continuously and the above estimate holds in Lebesgue spaces, i.e. 
	\begin{equation}
\|\partial^{\beta}E_{m,\alpha}(t)f\|_{L^{p_2}}\leq Ct^{-\frac{n}{2m}\big(\frac{1}{p_1}-\frac{1}{p_2}\big)-\frac{\alpha}{2m}-\frac{|\beta|}{2m}}\|f\|_{L^{p_1}}.
\end{equation}
	\end{enumerate}
\end{prop}

\begin{proof}
Denote by $1_{E}$ the characteristic function of the set $E\subset \mathbb{R}^{n}$. We may write  $$|x|^{-\alpha}=|x|^{-\alpha}1_{\{|x|>1\}}+|x|^{-\alpha}1_{\{|x|\leq 1\}}=f_1(x)+f_2(x)$$ and for $f\in L^{p_1}(\mathbb{R}^{n})$, $E_{m,\alpha}(t)f=E_{m,0}(t)(f_1f)+E_{m,0}(t)(f_2f)$ for any $t>0$. It is clear that $f_1\in L^{k_1}(\RR^{n})$, $f_2\in L^{k_2}(\RR^{n})$ whenever $k_1\alpha>n$ and $k_2\alpha<n$. We may then take $k_1=\frac{n}{\alpha}+\theta$ and $k_2=\frac{n}{\alpha}-\delta>0$ for some $\theta,\delta>0$. On the other hand, we have by H\"{o}lder's inequality $f_1f\in L^{l_1}(\RR^{n})$, $f_2f\in L^{l_2}(\RR^{n})$ provided $1\leq l_1,l_2\leq p_2$ $(l_1,l_2 \neq \infty)$ satisfy $\frac{1}{l_1}=\frac{1}{k_1}+\frac{1}{p_1}$, $\frac{1}{l_2}=\frac{1}{k_2}+\frac{1}{p_1}$. This shows that both $E_{m,0}(t)(f_1f)$ and $E_{m,0}(t)(f_2f)$ belongs to $C_0(\mathbb{R}^{n})$ and so is $E_{m,\alpha}(t)f$ for all $t>0$. This proves part \ref{part a}. To establish \ref{part b} note that
	by the continuous embedding $L^{p_2,1}(\mathbb{R}^{n})\hookrightarrow L^{p_2,q}(\mathbb{R}^{n})$, $1\leq q\leq \infty$, $p_2<\infty$ (\ref{eq:Lorentz-smoothing-effect}) easily follows from the following bound  
	$$\|\partial^{\beta}E_{m,\alpha}(t)f\|_{L^{p_2,1}}\leq Ct^{-\frac{n}{2m}\big(\frac{1}{p_1}-\frac{1}{p_2}\big)-\frac{\alpha}{2m}-\frac{|\beta|}{2m}}\|f\|_{L^{p_1,\infty}}\quad \text{for all}\quad t>0,\hspace{0.1cm}f\in L^{p_1,\infty}$$
	which we now prove. Let $1<k,l<\infty$ such that $1+\frac{1}{p_2}=\frac{1}{k}+\frac{1}{l}$. Using the generalized Young inequality (Lemma \ref{lem:Young ineq}), it follows that
	\begin{align*}
	\|\partial^{\beta}E_{m,\alpha}(t)f\|_{L^{p_2,1}}&=\|\partial^{\beta}e^{-t(-\Delta)^{m}}|\cdot|^{-\alpha}f\|_{L^{p_2,1}}\\
	&\leq  C\|\partial^{\beta}g_{m}(\cdot,t)\|_{L^{k,1}}\||\cdot|^{-\alpha}f\|_{L^{l,\infty}}\\
	&\leq C\|t^{-\frac{n}{2m}}\partial^{\beta}g(t^{-\frac{1}{2m}}\cdot)\|_{L^{k,1}}\||\cdot|^{-\alpha}\|_{L^{\frac{n}{\alpha},\infty}}\|f\|_{L^{p_1,\infty}}
	\end{align*}
	where in the last estimate, we have utilized the generalized H\"{o}lder's inequality (Lemma \ref{lem:Holder ineq}) with $\frac{1}{l}=\frac{1}{n/\alpha}+\frac{1}{p_1}<1$. Thus given that $g\in L^{k,1}(\RR^{n})$ in view of the pointwise bound (\ref{ineq:pointwise-on-g}), the scaling law of the Lorentz norm produces 
	\begin{align*}
	\|\partial^{\beta}E_{m,\alpha}(t)f\|_{L^{p_2,1}}&\leq Ct^{-\frac{n}{2m}}t^{-\frac{|\beta|}{2m}}t^{\frac{n}{2m k}}\|\partial^{\beta}g\|_{L^{k,1}}\|f\|_{L^{p_1,\infty}}\\
	&\leq Ct^{-\frac{n}{2m}(1-\frac{1}{k})}t^{-\frac{|\beta|}{2m}}\|f\|_{L^{p_1,\infty}}\\
	&\leq Ct^{-\frac{n}{2m}(\frac{1}{n/\alpha}+\frac{1}{p_1}-\frac{1}{p_2})}t^{-\frac{|\beta|}{2m}}\|f\|_{L^{p_1,\infty}}\\
	&\leq Ct^{-\frac{n}{2m}(\frac{1}{p_1}-\frac{1}{p_2})-\frac{\alpha}{2m}}t^{-\frac{|\beta|}{2m}}\|f\|_{L^{p_1,\infty}}.
	\end{align*}
	The remaining bit of the proof immediately flows from the continuous embeddings of $L^{p_2,1}(\RR^{n})$ into $L^{p_2,p_2}(\RR^{n})=L^{p_2}(\RR^{n})$ and $L^{p_1}(\RR^{n})\hookrightarrow L^{p_1,\infty}(\RR^{n})$, $p_1<\infty$. 
\end{proof}

\begin{remark} Observe that the estimate (\ref{eq:Lorentz-smoothing-effect}) in case $p_2=\infty$ remains valid but only under the condition $q=\infty$. By the possibility of choosing $\frac{1}{p_2}=\frac{1}{p_1}+\frac{\alpha}{n}$, we see that $p_1$ can be taken larger than $p_2$. This fact, of course is forbidden if $\alpha=0$. Notice that (\ref{ineq:condp-q}) may be rephrased as $\frac{1}{p_2}\leq \frac{\alpha}{n}+\frac{1}{p_1}<1$.   
\end{remark}

Let $F$ as in (\ref{eq:cond-on-F}). We introduce the nonlinear operator
\[
\mathcal{M}\phi(x)=\int_{0}^{\infty}E_{m,\alpha}(s)F(\phi(\cdot,s))ds
\]
whenever the integral on the right hand side makes sense for a suitable $\phi$. We show below that $\mathcal{M}$ is well-behaved between  Marcinkiewicz spaces.
\begin{lemma}\label{lem:nonlinearbound}
Let $p>1$ and $0<\alpha<2m<n$. Given $d,q>1$, assume that \[d>\frac{np}{n-\alpha}\hspace{0.2cm}\text{and}\hspace{0.2cm} \frac{1}{q}=\frac{p}{d}-\frac{2m-\alpha}{n}.\] There exists a constant $C>0$ such that 
\begin{equation}\label{eq:1st-bound-on-M}
\|\mathcal{M}\phi\|_{L^{q,\infty}}\leq C\sup_{t>0}\|\phi(t)\|^{p}_{L^{d,\infty}}, \quad \phi\in L^{\infty}([0,\infty);L^{d,\infty}(\RR^{n})).
\end{equation} 
In particular we can take $d=q_c$ in (\ref{eq:1st-bound-on-M}) and if $\phi\in L^{\infty}([0,\infty);L^{q_c,\infty}\cap L^{r,\infty}(\RR^{n}))$ with $r>1$ and $1<r'<\frac{n}{2m-\alpha}$, $\frac{1}{r}+\frac{1}{r'}=1$, then we have 
\begin{equation}\label{eq:snd-bound-on-M}
\|\mathcal{M}\phi\|_{L^{r,\infty}}\leq c\sup_{t>0}\|\phi(t)\|_{L^{r,\infty}}\sup_{t>0}\|\phi(t)\|^{p-1}_{L^{q_c,\infty}}.
\end{equation} 
\end{lemma}
\begin{proof}[Proof of Lemma \ref{lem:nonlinearbound}]
Denote by $\langle\cdot,\cdot\rangle$ the duality bracket between $L^{q,\infty}$ and $L^{q',1}$, $q>1$ where $q'$ is the conjugate exponent of $q$. Let $\psi\in L^{q',1}$ and $\phi\in L^{\infty}([0,\infty);L^{d,\infty}(\RR^{n}))$, $d>1$. We use Fubini's theorem and the generalized H\"{o}lder inequality to get 
\begin{align*}
\|\mathcal{M}\phi\|_{L^{q,\infty}}&=\sup_{\substack{\psi\in L^{q',1}\\
	\|\psi\|_{L^{q',1}}=1}}\big|\langle\mathcal{M}\phi,\psi\rangle\big|\\
&\leq \sup_{\substack{\psi\in L^{q',1}\\
		\|\psi\|_{L^{q',1}}=1}}\bigg| \int_{0}^{\infty}\int_{\RR^n}E_{m,\alpha}(s)F(\phi(x,s))\psi(x)dxds\bigg|\\
&\leq \sup_{\substack{\psi\in L^{q',1}\\
		\|\psi\|_{L^{q',1}}=1}}\bigg| \int_{0}^{\infty}\int_{\RR^n}E_{m,0}(s)|\cdot|^{-\alpha}F(\phi(x,s))\psi(x)dxds\bigg|\\	
&\leq \sup_{\substack{\psi\in L^{q',1}\\
		\|\psi\|_{L^{q',1}}=1}}\bigg(\int_{0}^{\infty}\|E_{m,0}(s)\psi\|_{L^{q_1,1}}\||\cdot|^{-\alpha}F(\phi(s))\|_{L^{q_1',\infty}}ds\bigg)\\
&\leq C\sup_{\substack{\psi\in L^{q',1}\\
		\|\psi\|_{L^{q',1}}=1}}\bigg(\int_{0}^{\infty}\|E_{m,0}(s)\psi\|_{L^{q_1,1}}\||\cdot|^{-\alpha}|\phi(s)|^{p}\|_{L^{q'_1,\infty}}ds\bigg)\\
&\leq C\sup_{\substack{\psi\in L^{q',1}\\
		\|\psi\|_{L^{q',1}}=1}}\bigg(\int_{0}^{\infty}\|E_{m,0}(s)\psi\|_{L^{q_1,1}}\||\cdot|^{-\alpha}|\|_{L^{\frac{n}{\alpha},\infty}}\|\phi(s)\|^{p}_{L^{d,\infty}}ds\bigg)\\
&\leq C\sup_{s>0}\|\phi(s)\|^p_{L^{d,\infty}}\bigg(\sup_{\substack{\psi\in L^{q',1}\\
		\|\psi\|_{L^{q',1}}=1}}\int_{0}^{\infty}\|E_{m,0}(s)\psi\|_{L^{q_1,1}}ds\bigg)
	\end{align*}
where $q_{1}$ is related to the other parameters according to $\frac{1}{q_1}= \frac{n-\alpha}{n}-\frac{p}{d}$. The fact that $q_1'>1$ yields the restriction $\frac{np}{n-\alpha}<d$. Put $\frac{1}{q}=\frac{p}{d}-\frac{2m-\alpha}{n}$, one easily verifies that $\frac{n}{2m}(\frac{1}{q'}-\frac{1}{q_1})-1=0$. Now since for all $1<p_1<p_2<\infty$, the following estimate holds (see e.g. Lemma 3.10 in \cite{Ferreira-Villa}) 
\[\int_{0}^{\infty}s^{\frac{n}{2m}(\frac{1}{p_1}-\frac{1}{p_2})-1}\|E_{m,0}(s)\phi\|_{L^{p_2,1}}ds\leq C\|\phi\|_{L^{p_1,1}}
\]
for all $\phi\in L^{p_1,1}(\RR^{n})$, it suffices to make the choices $p_1=q'$ and $p_2=q_1$ and we immediately find that $\displaystyle
\int_{0}^{\infty}\|E_{m,0}(s)\psi\|_{L^{q_1,1}}ds\leq c\|\psi\|_{L^{q',1}}$ from which we deduce the desired bound (\ref{eq:1st-bound-on-M}). We prove estimate (\ref{eq:snd-bound-on-M}) in a similar fashion. Indeed let $r>1$  and take $q_2>1$ with $\frac{1}{q_2}=\frac{1}{r}+\frac{2m-\alpha}{n}$. For every $t>0$, one has 
\[\|F(\phi(t))\|_{L^{q_2,\infty}}\leq C_1\||\phi(t)|^{p-1}\|_{L^{\frac{n}{2m-\alpha}}}\|\phi(t)\|_{L^{r,\infty}}=C_1\|\phi(t)\|^{p-1}_{L^{q_c},\infty}\|\phi(t)\|_{L^{r,\infty}}.
\] 
From $q_2^{-1}<1$, we obtain $r'<\frac{2m-\alpha}{n}$ and since $\frac{n}{2m}\big(\frac{1}{r'}+\frac{n-\alpha}{n}+\frac{1}{q_2}\big)=1$, we can simply pick $\psi\in L^{r',1}(\RR^{n})$ and mimic the above steps to reach (\ref{eq:snd-bound-on-M}). The proof of Lemma \ref{lem:nonlinearbound} is now complete.
\end{proof}

\section{Proofs of the main results}
\begin{proof}[Proof of Theorem \ref{thm:LWP-positive-alpha}]
	Let $\alpha<0$. Assume without loss of generality that $T\leq 1$ and consider for $M>0$ (to be chosen later), the space
	\[X_{T}=\{u\in L^{\infty}([0,T];C_{\Lambda}(\RR^{n})):\|u(t)\|_{\Lambda}\leq M, \hspace{0.1cm}t\in [0,T]\}
	\] 
	which when equipped with the metric $d(u,v)=\displaystyle\sup_{t\in [0,T]}\|u(t)-v(t)\|_{\Lambda}$ carries out the structure of a complete metric space. Next, consider the map $\mathcal{N}$ defined as
	\[\mathcal{N}u(t)=E_{m,0}(t)u_0+\int_{0}^{t}E_{m,\alpha}F(u(s))ds,\quad 0<t<T.
	\]
	Applying Lemma \ref{lem:exp-est-g}, we are able to estimate each of the terms in $\mathcal{N}$ as follows:
	\begin{align*}
	|E_{m,0}(t)u_0|&\leq \bigg|\int_{\RR^{n}}t^{-n/2m}g(t^{-1/2m}(x-y))u_0(y)dy\bigg|\\
	&\leq \bigg|\int_{\RR^{n}}t^{-n/2m}g(t^{-1/2m}(x-y))(1+|y|)^{-\frac{\alpha}{p-1}}\Lambda(y)u_0(y)dy\bigg|\\
	&\leq \|u_0\|_{\Lambda}\int_{\RR^{n}}|g(z)|(1+|x-t^{1/2m}z|)^{-\frac{\alpha}{p-1}}dy\\
	&\leq C(n,p,\alpha)\|u_0\|_{\Lambda}(1+|x|)^{-\frac{\alpha}{p-1}}
	\end{align*}
	where we have made the change of variables $x-y=t^{1/2m}z$. Also,
	\begin{align*}
	\bigg|\int_{0}^{t}E_{m,\alpha}F(u(s))ds\bigg|&= \bigg|\int_{0}^{t}\int_{\RR^{n}}(t-s)^{-n/2m}g\big((t-s)^{-1/2m}(x-y)\big)|y|^{\alpha}F(u(y,s))dy\bigg|\\
	&=\bigg| \int_{0}^{t}\int_{\RR^{n}}(t-s)^{-n/2m}g\big((t-s)^{-1/2m}(x-y)\big)|y|^{\alpha}\Lambda^{-1}(y)(\Lambda F(u(s)))(y)dy\bigg|\\
	&\leq\|u\|^{p}_{X_T} \int_{0}^{t}\int_{\RR^{n}}|g(z)|(1+|x-t^{1/2m}z|)^{-\frac{\alpha}{p-1}}dyds\\
	&\leq C_1(n,p,\alpha)T\|u\|^{p}_{X_T}\Lambda^{-1}(x)
	\end{align*}
	so that $\displaystyle \sup_{0\leq t\leq T}\|\mathcal{N}u\|_{X_{T}}\leq (C_1(n,p,\alpha)M^{p}T+C(n,p,\alpha)\|u_0\|_{\Lambda})$. Likewise, for any $u,v\in X_{T}$, we get 
	\[\|\mathcal{N}u-\mathcal{N}v\|_{X_{T}}\leq C_1(n,\alpha,p)M^{p-1}T\|u-v\|_{X_{T}}.
	\]
 For $T$ sufficiently small depending on $u_0$ we can achieve $C_1(n,p,\alpha)M^{p-1}T<1$ if $M>0$ is chosen small enough. Let  $\|u_0\|_{\Lambda}<A$ for $A>0$ chosen small such that $C_1(n,p,\alpha)M^{p}T+C(n,p,\alpha)A\leq M$ from which one deduces the self-mapping and contraction properties of the map $\mathcal{N}$. The Banach fixed point theorem yields the existence of exactly one mild solution $u$ of \eqref{eq:main-eq} in $X_{T}$. Given that $C_{\Lambda}(\RR^{n})\subset C_{b}(\RR^{n})$, we have that $u\in L^{\infty}([0,T];C_{b}(\RR^{n}))$ and by parabolic regularity theory, $u$ satisfies Eq. (\ref{eq:main-eq}) in the classical sense. 
\end{proof}

\begin{proof}[Proof of Theorem \ref{thm:GWP-Lp}]\label{subsec:proof-GWP}
Consider $u_0\in L^{q_c,\infty}(\RR^{n})$ such that $\|u_0\|_{L^{q_c,\infty}}<\varepsilon$ for some $\varepsilon>0$. Put $w(t)=E_{m,0}(t)u_0$ and define for $t>0$, the map 
\[Tu(t)=\int_{0}^{t}E_{m,\alpha}(t-s)F(u(s))ds.
\]
Let $p\in (\frac{n-\alpha}{n-2m},\infty)$. Apply Lemma \ref{lem:nonlinearbound} with the choice $d=q_c$ (i.e. $q=q_c$) in (\ref{eq:1st-bound-on-M})  which gives the estimate
\begin{equation}\label{self-map}[Tu]_{L^{q_c,\infty}}\leq K[u]^p_{L^{q_c,\infty}}.
\end{equation}
Moreover, given $u$ and $v$ in $L^{\infty}([0,\infty);L^{q_c,\infty}(\rn))$, we can proceed as in the proof of Lemma \ref{lem:nonlinearbound} to get 
\begin{equation}\label{contraction}
[T(u)-T(v)]_{L^{q_c,\infty}}\leq K[u-v]_{L^{q_c,\infty}}([u]^{p-1}_{L^{q_c,\infty}}+[v]^{p-1}_{L^{q_c,\infty}}).
\end{equation}
 Let $\vartheta=C\varepsilon$, where $C$ is the constant in the estimate $\|w\|_{L^{q_c,\infty}}\leq C\|u_0\|_{L^{q_c,\infty}}$. For $R=(2^{p}K)^{\frac{1}{1-p}}$ and by assuming that $\vartheta <R$, we find that $T$ is a self-mapping and a contraction on the closed ball $B_{2\vartheta}(0)=\{u\in L^{\infty}([0,\infty);L^{q_c,\infty}(\rn)):[u]_{L^{q_c,\infty}}\leq 2\vartheta\}$. Hence, $T$  has a unique fixed in $B_{2\vartheta}(0)$.  Next, we show the convergence of $u(t)$ to $u_0$ in the weak-$\star$ sense in $L^{q_c,\infty}(\RR^{n})$. Let $\varphi\in L^{q'_c,1}(\RR^{n})$, we have 
  \begin{align*}
  |\langle u(t)-u_0,\varphi\rangle|&=|\langle w-u_0,\varphi\rangle|+|\langle T(u),\varphi\rangle|\\
&= |\langle (E_{m,0}(t)\varphi-\varphi),u_0\rangle|+|\langle T(u),\varphi\rangle|\\
&=\|u_0\|_{L^{q_c,\infty}}\|E_{m,0}(t)\varphi-\varphi\|_{L^{q'_c,1}}+|\langle T(u),\varphi\rangle|.
  \end{align*}
Since $q'_c>1$, it is clear that $\|E_{m,0}(t)\varphi-\varphi\|_{L^{q'_c,1}}=\|g_{m}(\cdot,t)\ast\varphi-\varphi\|_{L^{q'_c,1}}\rightarrow 0$ as $t\rightarrow 0^{+}$. On the other hand, (\ref{eq:Lorentz-smoothing-effect}) of Proposition  \ref{prop:smoothing-effect}, applied with $(p_1,p_2)=(q_c/p,q_c)$ and $q=\infty$ permits us to write
\begin{align*}
|\langle T(u)(t),\varphi\rangle|&=\bigg|\int_{0}^{t}\int_{\RR^{n}}E_{m,\alpha}(s)\varphi(x)F(u(t-s))dxds\bigg|\\
&\leq \int_{0}^{t}\|E_{m,\alpha}(s)F(u(t-s))\|_{L^{q_c,\infty}}\|\varphi\|_{L^{q'_c,1}}ds\\
&\leq C\|\varphi\|_{L^{q'_c,1}}\sup_{t>0}\|u(t)\|^{p}_{L^{q_c,\infty}}\int_{0}^{t}s^{-(1-\frac{\alpha}{2m})}ds\\
&\leq C(2\vartheta)^{p}\|\varphi\|_{L^{q'_c,1}}t^{\alpha/2m}
\end{align*}    
which shows that $|\langle T(u)(t),\varphi\rangle|$ converges to $0$ when $t\rightarrow 0^{+}$ as desired. To prove the second part \ref{GWP-P2} of the theorem, note that the solution constructed above can be realized as the limit (in the space used above) of the successive approximations $(u_i)$ with
\[u_1(\cdot,t)=E_{m,0}(t)u_0,\quad u_{i+1}(t)=\int_{0}^{t}E_{m,\alpha}(t-s)F(u_i(s))ds+u_1,\hspace{0.1cm}i\in \NN_{\geq 2}.
\]
Take $\bar{\varepsilon}>0$, $\Bar{\varepsilon}\leq \varepsilon$ such that $2^{p}c\bar{\vartheta}^{p-1}<1$ where $c$ is the constant appearing in (\ref{eq:snd-bound-on-M}) and $\bar{\vartheta}=C\bar{\varepsilon}$. Let $u_0\in L^{q_c,\infty}\cap L^{r,\infty}(\RR^{n})$ with $\|u_0\|_{L^{q_c,\infty}}<\bar{\varepsilon}$. We claim that the above iteration sequence converges to a function $\bar{u}$ which uniquely solves $(\ref{eq:mild-solution})$. To justify the latter, we show that $(u_i)$ is a Cauchy sequence in the corresponding space. First notice that $u_i\in L^{\infty}([0,\infty);L^{r,\infty}(\rn))$ for each $i\in \NN$, in view of the assumptions imposed on $u_0$ and Lemma \ref{lem:nonlinearbound}. Furthermore, we invoke once again the latter Lemma which guarantees that
\[\|u_1\|_{L^{r,\infty}}<C'\Bar{\vartheta},\quad [u_{i+1}-u_i]_{L^{r,\infty}}\leq c2^{p}\bar{\vartheta}^{p-1}\sup_{t>0}\|(u_{i}-u_{i-1})(t)\|_{L^{r,\infty}}.
\] Iterating the second inequality leads to the estimate 
\[[u_{i+1}-u_i]_{L^{r,\infty}}\leq (c2^{p}\bar{\vartheta}^{p-1})^{i-2}\sup_{t>0}\|(u_{2}-u_{1})(t)\|_{L^{r,\infty}}.
\] 
Since $[u_{2}-u_1]_{L^{r,\infty}}<\infty$, we deduce that $u_{i+1}-u_i$ converges to $0$ as $i\rightarrow\infty$. Thus, $(u_i)$ is Cauchy and by uniqueness of the limit in the sense of distributions, we obtain the claim.\\
\textbf{Proof of Part \ref{GWP-P3}}.  Take $\alpha>0$ such that $\alpha+2m<n$ and assume that $p\in (\frac{n-\alpha}{n-2m},\frac{2m}{\alpha})$. For $q\in (q_c,\frac{np(p-1)}{2m-\alpha p})$, set $\delta=\frac{n}{2m}(1/q_c-1/q)$. Next, consider the Banach space $S_q^{\delta}$, defined as  
\[S_q^{\delta}=\{u\in C_{b}([0,\infty);L^{q_c,\infty}(\RR^{n})):\sup_{t>0}t^{\delta}\|u(t)\|_{L^{q,\infty}}<\infty\}\]  
and equipped with the norm $\|u\|_{S_q^{\delta}}=[u]_{L^{q_c,\infty}}+\sup_{t>0}t^{\delta}\|u(t)\|_{L^{q,\infty}}$.
Let $\varepsilon_q>0$ such that $\|u_0\|_{L^{q_c,\infty}}<\varepsilon_q$. By the smoothing effect, we have 
\begin{align*}
\|E_{m,0}(t)u_0\|_{L^{q,\infty}}\leq Ct^{-\frac{n}{2m}(\frac{1}{q_c}-\frac{1}{q})}\|u_0\|_{L^{q_c,\infty}}\leq C't^{-\delta}\varepsilon_q.
\end{align*}
Consider the operator $T$ previously defined, it satisfies 
\begin{align*}
\|T(u)(t)\|_{L^{q,\infty}}&\leq \int_{0}^{t}\|E_{m,\alpha}(t-s)F(u(s))\|_{L^{q,\infty}}ds\\
&\leq c\int_{0}^{t}(t-s)^{-\frac{n}{2m}(\frac{p-1}{q})-\frac{\alpha}{2m}}s^{-\delta p}s^{\delta p}\|u(s)\|^{p}_{L^{q,\infty}}ds\\
&\leq c(\sup_{t>0}t^{\delta}\|u(t)\|)^{p}_{L^{q,\infty}}\int_{0}^{t}(t-s)^{-\frac{n}{2m}(\frac{p-1}{q})-\frac{\alpha}{2m}}s^{-p\delta}ds\\
&\leq c(\sup_{t>0}t^{\delta}\|u(t)\|)^{p}_{L^{q,\infty}}t^{-\delta p-\frac{n(p-1)}{2mq}-\frac{\alpha}{2m}+1}\int_{0}^{1}s^{-\delta p}(1-s)^{-\frac{n(p-1)}{2mq}-\frac{\alpha}{2m}}ds
\end{align*}  
where we have used in the second inequality Proposition \ref{prop:smoothing-effect} together with $\frac{1}{q}<\frac{n-\alpha}{np}$ i.e. $\frac{p}{q}+\frac{\alpha}{n}<1$ which holds since $q>q_c$ and $p>\frac{n-\alpha}{n-2m}$. Since $\delta p+\frac{n(p-1)}{2mq}+\frac{\alpha}{2m}+1=\delta$, $\delta p<1$ and $1-\frac{n(p-1)}{2mq}-\frac{\alpha}{2m}>0$, it easily follows that 
\[\|Tu(t)\|_{L^{q,\infty}}\leq C(\sup_{t>0}t^{\delta}\|u(t)\|_{L^{q,\infty}})^{p}t^{-\delta}\mathcal{B}(1-\delta p,1-\frac{n(p-1)}{2mq}-\frac{\alpha}{2m})\] where  $\mathcal{B}(s,\tau)=\displaystyle
\int_{0}^{1}\sigma^{s-1}(1-\sigma)^{\tau-1}d\sigma=\frac{\Gamma(s)\Gamma(\tau)}{\Gamma(s+\tau)},$ $s,\tau>0$ and $\Gamma$ is the Gamma function.
Following the lines of the above argument combined with the proof of Part \ref{GWP-P1} leads to
\[\|Tu\|_{S^{\delta}_q}\leq C''\|u\|^{p}_{S^{\delta}_q},\quad\|Tu-Tv\|_{S^{\delta}_q}\leq C''\|u-v\|_{S^{\delta}_q}\|u\|^{p-1}_{S^{\delta}_q},
\] 
for all $u,v\in S^{\delta}_q$. Now, pick $R=(2^{p}C'')^{\frac{1}{1-p}}$ and assume $\varepsilon_q<R/2C'$; we then have $\|E_{m,0}(t)u_0\|_{S^{\delta}_q}<\vartheta_q$, $\vartheta_q=C'\varepsilon_q$. Applying the Banach fixed point theorem,  the conclusion is straightforward. 
\end{proof}

\begin{proof}[Proof of the uniqueness criterion (Theorem \ref{thm:uniqueness})] Let $p\in(\frac{n-\alpha}{n-2m},\infty)$ where $0<\alpha<2m<n$ and consider $u,v$ two  global mild solutions  in $C_{b}([0,\infty);L^{q_c}(\RR^{n}))$ initially given by $u_0\in L^{q_c}(\RR^{n})$. We show that $u(t)=v(t)$ for $t\in [0,\infty)$. Pick $T>0$ small. We start by showing that $u$ and $v$ agree on $[0,T]$. For this end, set $w(t)=u(t)-v(t)$ and use (\ref{eq:cond-on-F}) to write 
\begin{align*}
\|w(t)\|_{L^{q_c,\infty}}&=\bigg\|\int_{0}^{t}E_{m,\alpha}(t-s)(F(u(s))-F(v(s))ds\bigg\|_{L^{q_c,\infty}}\\
&\leq C\bigg\|\int_{0}^{t}E_{m,\alpha}(t-s)|w(s)||u(s)-E_{m,0}(s)u_0|^{p-1}ds\bigg\|_{L^{q_c,\infty}}+\\
&\hspace{3.4cm}+C\bigg\|\int_{0}^{t}E_{m,\alpha}(t-s)|w(s)||v(s)-E_{m,0}(s)u_0|^{p-1}\big)ds\bigg\|_{L^{q_c,\infty}}+\\
&\hspace{5cm}2C\bigg\|\int_{0}^{t}E_{m,\alpha}(t-s)|w(s)||E_{m,0}(s)u_0|^{p-1}ds\bigg\|_{L^{q_c,\infty}}\\
&= \mathcal{I}+\mathcal{J}+\mathcal{K}.
\end{align*}
Imitating the proof of Lemma \ref{lem:nonlinearbound} with $(q,d)=(q_c,\frac{q_c}{p})$, under the condition $\frac{n-\alpha}{n-2m}<p<\infty$, $0<\alpha<2m<n$ gives
\begin{align}\label{eq:bound-on-I}\nonumber \mathcal{I}&\leq C\sup_{t\in [0,T]}\||w(t)||u(t)-E_{m,0}(t)u_0|^{p-1}\|_{L^{q_c/p,\infty}}\\
\nonumber &\leq C\sup_{0\leq t\leq T}\|w(t)\|_{L^{q_c,\infty}}\sup_{0\leq t\leq T}\|u(t)-E_{m,0}(t)u_0\|^{p-1}_{L^{q_c,\infty}}\\
&\leq C\sup_{0\leq t\leq T}\|w(t)\|_{L^{q_c,\infty}}\sup_{0\leq t\leq T}\|u(t)-E_{m,0}(t)u_0\|^{p-1}_{L^{q_c}}
\end{align}
where we have utilized the generalized H\"{o}lder's inequality with $p/q_c=(p-1)/q_c+1/q_c$. Similarly, it holds that
\begin{equation}\label{eq:bound-on-J}
\mathcal{J}\leq C\sup_{0\leq t\leq T}\|w(t)\|_{L^{q_c,\infty}}\sup_{0\leq t\leq T}\|v(t)-E_{m,0}(t)u_0\|^{p-1}_{L^{q_c}}.
\end{equation}
 Now take $r>q_c$ and put $\delta=\frac{n}{2m}(\frac{1}{q_c}-\frac{1}{r})$, $r_1^{-1}=\frac{p-1}{r}+\frac{1}{q_c}$. One easily verifies (taking the range of $p$ into consideration) the relations
 \begin{equation}\label{eq:relation-on-exponents}
 \frac{p-1}{r_1}+\frac{\alpha}{n}<1,\hspace{0.1cm} \delta(p-1)<1,\hspace{0.1cm} 1-\frac{n(p-1)-\alpha r}{2mr}>0, \hspace{0.1cm}\delta^{\star}=1-\frac{n}{2m}\big(\frac{1}{r_1}-\frac{1}{q_c}\big)-\frac{\alpha}{2m}-\delta(p-1)=0.
 \end{equation}  
 Finally we bound $\mathcal{K}$ by applying Proposition \ref{prop:smoothing-effect} with $q=\infty$ and H\"{o}lder inequality to get 
 \begin{align}\label{eq:bound-on-K}
 \nonumber\mathcal{K}&\leq c\int_{0}^{t}(t-s)^{-\frac{n(p-1)}{2mr}-\frac{\alpha}{2m}}\||w(t)||E_{m,0}(s)u_0|^{p-1}\|_{L^{r_1,\infty}}\\
 \nonumber&\leq c\sup_{0\leq t\leq T}\|w(t)\|_{L^{q_c,\infty}}(\sup_{0\leq t\leq T}t^{\delta}\|E_{m,0}(t)u_0\|_{L^{r}})^{p-1}\int_{0}^{t}(t-s)^{-\frac{n(p-1)}{2mr}-\frac{\alpha}{2m}}s^{-\delta(p-1)}ds\\
 \nonumber&\leq c\sup_{ t\in [0,T]}\|w(t)\|_{L^{q_c,\infty}}(\sup_{ t\in [0,T]}t^{\delta}\|E_{m,0}(t)u_0\|_{L^{r}})^{p-1}t^{\delta^{\star}}\int_{0}^{1}s^{-\delta(p-1)}(1-s)^{-\frac{n(p-1)}{2mr}-\frac{\alpha}{2m}}ds\\
  \nonumber&\leq c\sup_{ t\in [0,T]}\|w(t)\|_{L^{q_c,\infty}}(\sup_{ t\in [0,T]}t^{\delta}\|E_{m,0}(t)u_0\|_{L^{r}})^{p-1}\mathcal{B}\big(\delta(p-1),1-\frac{n(p-1)}{2rm}-\frac{\alpha}{2m}\big)\\
 &\leq c\sup_{ t\in [0,T]}\|w(t)\|_{L^{q_c,\infty}}(\sup_{ t\in [0,T]}t^{\delta}\|E_{m,0}(t)u_0\|_{L^{r}})^{p-1}
 \end{align}
where we have made use of (\ref{eq:relation-on-exponents}). Thus, from (\ref{eq:bound-on-I}), (\ref{eq:bound-on-J}) and (\ref{eq:bound-on-K}), we arrive at 
\begin{equation}\label{eq:est-lambda}
\sup_{ t\in [0,T]}\|w(t)\|_{L^{q_c,\infty}}\leq C_1\sup_{ t\in [0,T]}\|w(t)\|_{L^{q_c,\infty}}\lambda(T),
\end{equation}
with $\lambda(T)=\displaystyle\sup_{t\in[0, T]}\|u(t)-E_{m,0}(t)u_0\|_{L^{q_c}}+\displaystyle\sup_{t\in [0,T]}\|v(t)-E_{m,0}(t)u_0\|_{L^{q_c}}+(\sup_{t\in [0,T]}t^{\delta}\|E_{m,0}(t)u_0\|_{L^{r}})^{p-1}.$
It has already been established in the proof of Theorem \ref{thm:GWP-Lp} that if $u_0\in L^{q_c}(\RR^{n})$, then it holds that $\displaystyle\limsup_{t\rightarrow 0^{+}}\|u(t)-E_{m,0}(t)u_0\|_{L^{q_c}}= 0=\displaystyle\limsup_{t\rightarrow 0^{+}}\|v(t)-E_{m,0}(t)u_0\|_{L^{q_c}}$. However, it can be shown in parallel that $\displaystyle\limsup_{t\rightarrow 0^{+}}t^{\delta}\|E_{m,0}(t)u_0\|_{L^{r}}=0.$ Indeed, via a density argument, there exists a sequence $(u_{0N})_{N}\in L^{q_c}\cap L^{r}(\RR^{n})$ satisfying $u_{0N}\rightarrow u_0$, ($N\rightarrow \infty$) in $L^{q_c}(\RR^{n})$ and the estimate  $t^{\delta}\|E_{m,0}(t)u_{0N}\|_{L^{r}}\leq Ct^{\delta}\|u_{0N}\|_{L^{r}}$ so that $\displaystyle\lim_{t\rightarrow 0^{+}}t^{\delta}\|E_{m,0}(t)u_{0N}\|_{L^{r}}=0$. Moreover, we have
\begin{align*}
t^{\delta}\|E_{m,0}(t)u_{0}\|_{L^{r}}&\leq t^{\delta}\|E_{m,0}(t)(u_0-u_{0N})\|_{L^{r}}+t^{\delta}\|E_{m,0}(t)u_{0N}\|_{L^{r}}\\
&\leq C\|u_0-u_{0N}\|_{L^{q_c}}+t^{\delta}\|E_{m,0}(t)u_{0N}\|_{L^{r}}
\end{align*}
and by passing to the limit simultaneously as $t\rightarrow 0^{+}$ and $N\rightarrow \infty$, we obtain the claim.
 On the other hand, we can prove (see proof of Theorem \ref{thm:GWP-Lp}) that  $\displaystyle\lim_{t\rightarrow 0^{+}}\|E_{m,0}(t)u_0-u_0\|_{L^{q_c}}=0$. 
Hence, if one takes $T>0$ sufficiently small we can achieve $C_1\lambda(T)<\mu$ with $\mu\lll 1$ so that because of (\ref{eq:est-lambda}), we have $\|w(t)\|_{L^{q_c,\infty}}=0$ for every $t\leq T$. The rest of the proof goes as follows. Let $T_1=2T$ so that $T_1-T$ is small and consider the functions $\tilde{u}(x,t)=u(x,t+T)$, $\tilde{v}(x,t)=v(x,t+T)$, $(x,t)\in \RR^{n}\times [0,T]$ which verifies 
\[\tilde{u}(x,t)=E_{m,0}(t)u(T)+\int_{0}^{t}E_{m,\alpha}(t-\tau)F(\tilde{u}(\tau)d\tau,\hspace{0.1cm} \tilde{v}(x,t)=E_{m,0}(t)v(T)+\int_{0}^{t}E_{m,\alpha}(t-\tau)F(\tilde{v}(\tau)d\tau,\]
respectively. Since $u(T)=v(T)$, we can proceed as before to estimate $\widetilde{w}=\tilde{u}-\tilde{v}$ and obtain an analogue of (\ref{eq:est-lambda}) for $\widetilde{w}$. This will yield $\widetilde{w}=0$ on $[0,T]$, i.e. $u(t)=v(t)$, for all $t\in [T,2T]$. Iterating the preceding argument produces the desired conclusion on $(0,\infty)$.   
\end{proof}
\begin{proof}[Proof of Theorem \ref{thm:GWP-decay}]
Assume $\alpha\in (0,2m)$ with $(2m+\alpha)<n$, take $\ell\in (q_c,\infty)$ such that \begin{equation*}
\dfrac{1}{q_c}-\dfrac{2m}{n}<\dfrac{1}{\ell}<\dfrac{n-\alpha}{np} \hspace{0.2cm}\mbox{and set} \hspace{0.2cm}\delta=\dfrac{n}{2m}(1/q_c-1/\ell).
\end{equation*}
 We easily verify that $\delta=\frac{n(p-1)}{2m\ell}+p\delta+\frac{\alpha}{2m}-1$. For $K>0$, introduce the complete metric space $$E_K=\{u\in L^{\infty}((0,\infty);L^{\ell}(\RR^{n})):\sup_{t>0}t^{\delta}\|u(t)\|_{L^{\ell}}\leq K\}$$ endowed with the metric $d(u,v)=\displaystyle \sup_{t>0}t^{\delta}\|u(t)-v(t)\|_{L^{\ell}}$; $u,v\in E_K$ and define 
\[\mathcal{P}u(t)=E_{m,0}(t)u_0+\int_{0}^{t}E_{m,\alpha}(t-\tau)F(u(\tau))d\tau.
\]
Remark that \eqref{eq:Lorentz-smoothing-effect} holds for $\alpha=0$ so that with the choices $(p_1,p_2)=(q_c,\ell)$ and $q=\ell$, we get 
\begin{equation}\label{Eu0}
\|E_{m,0}(t)u_0\|_{L^{\ell}}\leq C_0t^{-\delta}\bigg\||\cdot|^{-\frac{2m-\alpha}{p-1}}\bigg\|_{L^{q_c,\infty}}\leq C_0t^{-\delta}
\end{equation} for every $t>0$ since $|u_0(x)|\leq c_0|x|^{-\frac{2m-\alpha}{p-1}}$. Since $\ell>np/(n-\alpha)$, Proposition \ref{prop:smoothing-effect}  can be used to derive the following bounds 
\begin{align}\label{Pself}
\nonumber\bigg\|\int_{0}^{t}E_{m,\alpha}(t-\tau)F(u(\tau))d\tau\bigg\|_{L^{\ell}}&\leq C\int_0^{t}\|E_{m,\alpha}(t-\tau)F(u(\tau))\|_{L^{\ell}}\\ 
\nonumber&\leq C\int_0^t (t-\tau)^{-\frac{n(p-1)}{2m\ell}-\frac{\alpha}{2m}}\|u(\tau)\|^p_{L^{\ell}}d\tau\\
\nonumber&\leq C(\sup_{t>0}t^{\delta}\|u(t)\|_{L^{\ell}})^{p}\int_0^t (t-\tau)^{-\frac{n(p-1)}{2m\ell}-\frac{\alpha}{2m}}\tau^{-p\delta}d\tau\\
\nonumber&\leq C(\sup_{t>0}t^{\delta}\|u(t)\|_{L^{\ell}})^{p}t^{-\delta}\int_0^1 (1-\tau)^{-\frac{n(p-1)}{2m\ell}-\frac{\alpha}{2m}}\tau^{-p\delta}d\tau\\
&\leq CK^{p}t^{-\delta}
\end{align}  and likewise
\begin{align}\label{Pcont}\bigg\|\int_{0}^{t}E_{m,\alpha}(t-\tau)[F(u(\tau))-F(v(\tau))]d\tau\bigg\|_{L^{r}}\leq CK^{p-1}t^{-\delta}\sup_{t>0}(t^{\delta}\|u(t)-v(t\|_{L^{r}}).
\end{align}
If $K$ is chosen small,  then the conditions $CK^{p-1}<1$ and $CK^{p}+C_0\leq K$ can be achieved since $C_0$ depends on $c_0$ which is sufficiently small by hypothesis. From \eqref{Eu0}, \eqref{Pself} and \eqref{Pcont}, we conclude that $\mathcal{P}$ is a contractive self-mapping onto $E_K$. Thus, $\mathcal{P}$ has a unique fixed point $u$ in $E_{K}$ solution of \eqref{eq:mild-solution}. 
\end{proof}

\begin{proof}[Proof of Theorem \ref{thm:nonexistence}] Here, we suppose $m\in \mathbb{N}$ and $2m<n$.
Let $u_0\in C_0(\rn)$ and assume that $u\in C([0,\infty);C_0(\rn))$ is a mild solution of Eq. \eqref{eq:main-eq} such that $|x|^{-\alpha}|u|^p\in C_0(\rn\setminus\{0\})$. We first show that $u$ is a weak solution. Take $\phi\in C^1([0,\infty);C^{2m}(\rn))$ with $\mbox{supp} \phi \subset [0,\infty)\times \rn$. Multiply Eq. \eqref{eq:mild-solution} by $\phi$ and integrate in space to get
\begin{align*}
\int_{\rn}u(x,t)\phi(t,x)dx=\int_{\rn}E_{m,0}(t)u_0(x)\phi(t,x)dx+\int_{\rn}\int_0^tE_{m,0}(t-\tau)|\cdot|^{-\alpha}|u|^p(\tau)d\tau \phi(t,x)dx.    
\end{align*}
Since $u_0$ is continuous and bounded, $E_{m,0}(t)$ is smooth and we can differentiate under the integral sign. Moreover, using the semigroup properties of the polyharmonic heat kernel and integration by parts we find that 
\begin{align*}
\partial_t\bigg(\int_{\rn}u(x,t)\phi(t,x)dx\bigg)&=\int_{\rn}\partial_t (E_{m,0}(t)u_0(x)\phi(t,x))dx+\\
&\hspace{1cm}\int_{\rn}\partial_t\bigg(\int_0^tE_{m,0}(t-\tau)|\cdot|^{-\alpha}|u|^p(\tau)d\tau \phi(t,x)\bigg)dx\\
&=A_1+A_2.
\end{align*}
We have
\begin{align*}
A_1&=\int_{\rn}[-(-\Delta)^{m}E_{m,0}(t)u_0(x)\phi(t,x)+E_{m,0}(t)u_0(x)\partial_t\phi(t,x)]dx\\
&=\int_{\rn}E_{m,0}(t)u_0(x)[-(-\Delta)^{m}\phi(t,x)+\partial_t\phi(t,x)]dx
\end{align*}
and
\begin{align*}
A_2&=\int_{\rn}|x|^{-\alpha}|u|^p\phi(x,t)dx+\int_{\rn}\int_0^t-(-\Delta)^m(E_{m,0}(t-\tau)|\cdot|^{-\alpha}|u|^p(\tau)d\tau)\phi(t,x) dx+\\
&\hspace{5cm}\int_{\rn}\bigg(\int_0^tE_{m,0}(t-\tau)|\cdot|^{-\alpha}|u|^p(\tau)d\tau\bigg)\partial_t\phi(t,x)dx\\
&=\int_{\rn}|x|^{-\alpha}|u|^p\phi(t,x)dx+\int_{\rn}\int_0^tE_{m,0}(t-\tau)|\cdot|^{-\alpha}|u|^p(\tau)d\tau[-(-\Delta)^m\phi(t,x)+\partial_t\phi(t,x)]dx
\end{align*}
so that 
\begin{align*}
A_1+A_2&=\int_{\rn}\bigg(\int_0^tE_{m,0}(t-\tau)|\cdot|^{-\alpha}|u|^p(\tau)d\tau+E_{m,0}(t)u_0(x)\bigg)[-(-\Delta)^{m}\phi(x,t)+\partial_t\phi(t,x)]dx+\\
&\hspace{8.5cm}\int_{\rn}|x|^{-\alpha}|u|^p\phi(t,x)dx\\
&=\int_{\rn}u(t,x)[-(-\Delta)^{m}\phi(x,t)+\partial_t\phi(t,x)]dx+\int_{\rn}|x|^{-\alpha}|u|^p\phi(t,x)dx.
\end{align*}
Integrating the above equality over $(0,\infty)$, we arrive at
\begin{align*}
\int_0^{\infty}\int_{\rn}u(t,x)[(-\Delta)^{m}\phi(x,t)-\partial_t\phi(t,x)]dxdt&=\int^{\infty}_0\int_{\rn}|x|^{-\alpha}|u|^p\phi(t,x)dxdt+\\
&\hspace{4cm}\int_{\rn}u(x,0)\phi(0,x)dx.    
\end{align*}
Next, we proceed via the rescaled test function method introduced in \cite{Miti-Poho}.
Consider the function $\psi\in C^{\infty}_{0}([0,\infty))$ such that
\[0\leq \psi\leq 1, \quad \psi(\tau)=\begin{cases}
1 \hspace{0.2cm}\text{if}\hspace{0.2cm}\tau\in [0,1]\\
0 \hspace{0.2cm}\text{if}\hspace{0.2cm}\tau\in [2,\infty).
\end{cases}
\]	
Suppose $u_0\in C_0(\rn)\cap L^{1}(\RR^{n})$ with $\displaystyle\int_{\RR^{n}}u_0(x)dx>0$. For $p\in (1,p_{F})$ and $R>0$, introduce the function $\psi_R(x,t)=\psi\big(\frac{t}{R}\big)^{\frac{p}{p-1}}\psi\big(\frac{|x|^{2m}}{R}\big)^{\frac{2mp}{p-1}}$. By way of contradiction, assume that $u$ is a global mild solution of (\ref{eq:main-eq}), then $u$ is a weak solution as shown above and we have 
\begin{align*}
\int_{\RR^{n}_{+}}|x|^{-\alpha}|u|^{p}\psi_{R}dxdt+\int_{\RR^{n}}u_{0}\psi_{R}(x,0)dx=\int_{\RR^{n}_{+}}u(-\partial_t\psi_R+(-\Delta)^{m}\psi_{R})dxdt.
\end{align*} 
By virtue of Young's inequality the integral $I_{1}=-\displaystyle
\int_{\RR^{n}_{+}}u\partial_t\psi_Rdxdt$ can be estimated as follows: 
\begin{align*}
I_{1}&\leq \frac{1}{2}\int_{\RR^{n}_{+}}|x|^{-\alpha}|u|^{p}\psi_{R}dxdt+C \int_{\RR^{n}_{+}}|x|^{\frac{\alpha}{p-1}}|\partial_t\psi^{\frac{p}{p-1}}(t/R)|^{p'}\psi\big(R^{-1}|x|^{2m}\big)^{\frac{2mp}{p-1}}\psi(t/R)^{-\frac{p}{(p-1)^{2}}}dxdt\\
&\leq \frac{1}{2}\int_{\RR^{n}_{+}}|x|^{-\alpha}|u|^{p}\psi_{R}dxdt+C_pR^{-p'} \int_{\RR^{n}_{+}}|x|^{\frac{\alpha}{p-1}}|\psi'(t/R)|^{p'}\psi\big(R^{-1}|x|^{2m}\big)^{\frac{2mp}{p-1}}dxdt\\
&\leq \frac{1}{2}\int_{\RR^{n}_{+}}|x|^{-\alpha}|u|^{p}\psi_{R}dxdt+C_pR^{-p'}R\bigg(\int_{0}^{2}|\psi'(s)|^{p'}ds\bigg)\bigg(\int_{\RR^{n}}|x|^{\frac{\alpha}{p-1}}\psi(R^{-1}|x|^{2m})^{\frac{2mp}{p-1}}dx\bigg)\\
&\leq \frac{1}{2}\int_{\RR^{n}_{+}}|x|^{-\alpha}|u|^{p}\psi_{R}+C_pR^{-p'+1+\frac{\alpha}{2m(p-1)}+\frac{n}{2m}}\bigg(\int_{0}^{2}|\psi'|^{p'}(s)ds\bigg)\int_{\{|z|< 2^{\frac{1}{2m}}\}}|z|^{\frac{\alpha}{p-1}}\psi(|z|^{2m})^{\frac{2mp}{p-1}}dz\\
&\leq \frac{1}{2}\int_{\RR^{n}_{+}}|x|^{-\alpha}|u|^{p}\psi_{R}dxdt+CR^{-p'+1+\frac{\alpha}{2m(p-1)}+\frac{n}{2m}}.
\end{align*}
In the first estimate, the fact that $\partial_t \psi^{p'}=(p'-1)\psi'\psi^{p'-1}$ contributes to eliminate the term with the negative exponent.
Note that in the estimate before the last we made the change of variables $s=t/R$ and $z=R^{-\frac{1}{2m}}x$, respectively. Moving on,
let $\psi_{2}(x)=\psi(\frac{|x|^{2m}}{R})^{\frac{2mp}{p-1}}$. It can be verified that for all $R>0$,  $|(-\Delta)^{m}\psi_{2}|\leq CR^{-1}\psi^{\frac{2m}{p-1}}(R^{-1}\cdot)$ and we bound the integral $I_{2}=\displaystyle\int_{\RR^{n}_{+}}u(-\Delta)^{m}\psi_{R}dxdt$ in the same way as above writing
\begin{align*}
I_{2}&\leq \bigg|\int_{\RR^{n}_{+}}u(-\Delta)^{m}\psi_{R}dxdt\bigg|\\
&\leq \frac{1}{2}\int_{\RR^{n}_{+}}|x|^{-\alpha}|u|^{p}\psi_{R}dxdt+C\int_{\RR^{n}_{+}}|x|^{\frac{\alpha}{p-1}}\psi_{R}^{-\frac{1}{p-1}}|(-\Delta)^{m}\psi_{R}(R^{-1}|x|^{2m})|^{p'}dxdt\\
&\leq \frac{1}{2}\int_{\RR^{n}_{+}}|x|^{-\alpha}|u|^{p}\psi_{R}dxdt+C\int_{\RR^{n}_{+}}|x|^{\frac{\alpha}{p-1}}\psi(t/R)^{\frac{p}{p-1}}\psi(|x|^{2m}/R)^{-\frac{2mp}{(p-1)^{2}}}|(-\Delta)^{m}\psi_{2}|^{p'}dxdt\\
&\leq \frac{1}{2}\int_{\RR^{n}_{+}}|x|^{-\alpha}|u|^{p}\psi_{R}dxdt+CR^{-p'+1}R^{\frac{\alpha}{2m(p-1)}}R^{\frac{n}{2m}}\bigg(\int_{0}^{2}\psi^{\frac{p}{p-1}}(s)ds\bigg)\int_{\{ |z|< 2^{1/2m}\}}|z|^{\frac{\alpha}{p-1}}dz\\
&\leq \frac{1}{2}\int_{\RR^{n}_{+}}|x|^{-\alpha}|u|^{p}\psi_{R}dxdt+CR^{-p'+1}R^{\frac{\alpha}{2m(p-1)}}R^{\frac{n}{2m}}.
\end{align*}
Summarizing, we find that 
\begin{equation}\label{eq:contracdiction-ineq}
\int_{\RR^{n}}u_0(x)\psi_{R}(x,0)dx\leq CR^{\frac{n}{2m}+\frac{\alpha}{2m(p-1)}-\frac{1}{p-1}}. 
\end{equation}
On the other hand, applying the Dominated Convergence Theorem, we have that 
\[\int_{\RR^{n}}u_0(x)\psi_{R}(x,0)dx\rightarrow \int_{\RR^{n}}u_0(x)dx\quad\text{as}\quad R\rightarrow \infty.
\]  
Hence, for $R>R_0$, $R_0$ sufficiently large, we have $\displaystyle \int_{\RR^{n}}u_0(x)\psi(R^{-1}|x|^{2m})dx>\frac{1}{4}\int_{\RR^{n}}u_0(x)dx$ so that by passing to the limit in (\ref{eq:contracdiction-ineq}) as $R\rightarrow \infty$, we arrive at $\displaystyle \int_{\RR^{n}}u_0(x)dx\leq 0$ since $1<p<p_{F}$. We have reached a contradiction in view of the assumption imposed on $u_0$. 
\end{proof}
\section{Long-time asymptotic behavior in the critical space and qualitative properties of solutions}
In this section, our goal is to investigate certain intrinsic properties of globally defined solutions obtained in Theorem \ref{thm:GWP-Lp}. They include self-similarity, radial symmetry, positivity  and large time behavior. 
\subsection{Self-similar solutions and further features}
Assume that $u$ defined on $\RR^{n}\times (0,\infty)$ is a global solution of (\ref{eq:main-eq}). For any $\sigma>0$, one easily verifies that the rescaled function 
\begin{equation}\label{eq:self-similarity}
u_{\sigma}(x,t)=\sigma^{\frac{2m-\alpha}{p-1}}u(\sigma x,\sigma^{2m}t)
\end{equation} 
formally solves (\ref{eq:main-eq}) provided $F(u(\sigma x,\sigma^{2m}t))=\sigma^{-\frac{(2m-\alpha)p}{p-1}}F(u_{\sigma})$ (this holds for instance if $F(u)=C|u|^{p}$ or $F(u)=C|u|^{p-1}u$ for a non-zero constant $C\in \RR$). By taking the limit as $t\rightarrow 0^{+}$, one sees that such solutions if they exist, should arise from initial data $u_{\sigma,0}=\sigma^{\frac{2m-\alpha}{p-1}}u_0(\sigma \cdot)$, homogeneous functions of degree $-\frac{2m-\alpha}{p-1}$. Thus, in the spirit of Theorem \ref{thm:GWP-Lp}, a solution which is initially given by $f(x)=\varepsilon_0|x|^{-\frac{2m-\alpha}{p-1}}$ for small enough $\varepsilon_0>0$ is globally defined and satisfies $u=u_{\sigma}$. Given that $f\in L^{q_c,\infty}(\RR^{n})$, this observation can be formulated in a more general context. 
\begin{theorem}[Self-similarity]\label{thm:self-similarity}
Let $\alpha>0$ and $p>1$ fall under the scope of Theorem \ref{thm:GWP-Lp}. Suppose $F$ is as above and satisfies \eqref{eq:cond-on-F}. Further, let $u_0\in L^{q_c,\infty}(\RR^{n})$ homogeneous of degree $-\frac{2m-\alpha}{p-1}$ with small norm. There exists a unique solution $u$ of Eq. \eqref{eq:mild-solution} verifying $u=u_{\sigma}$ a.e. on $\RR^{n}\times [0,\infty)$.
Moreover, if  $\|u_0\|_{L^{q_c,\infty}}<\varepsilon_q$ with $q>1$ and $\varepsilon_q>0$ as in part \ref{GWP-P3} in Theorem \ref{thm:GWP-Lp}, then the corresponding solution in $S^{\delta}_q$ is self-similar, i.e. $u=u_{\sigma}$.
\end{theorem}
\begin{corollary}\label{cor:trivial-self-similar}
There is no nontrivial self-similar solution of Eq. \eqref{eq:mild-solution} in $C_{b}([0,\infty);L^{q_c}(\RR^{n}))$ which is initially small in $L^{q_c}(\RR^{n})$.  
\end{corollary}
\begin{proof}[Proof of Theorem \ref{thm:self-similarity}]
Let $u\in C_{b}([0,\infty);L^{q_c,\infty}(\RR^{n}))$ be a global solution of Eq. \eqref{eq:mild-solution}. For any $\sigma>0$, we claim that $u_{\sigma}$ with $u_{\sigma}(\cdot,0)=u_{\sigma,0}$  also satisfies Eq. \eqref{eq:mild-solution}. Indeed, for all $t>0$, using $g_m(\cdot,t)=t^{n/2m}g(|\cdot|/t^{1/2m})$, we have
\begin{align*}
(E_{m,0}(t)u_0)_{\sigma}(x)&=\sigma^{\frac{2m-\alpha}{p-1}}\int_{\RR^{n}}g_{m}(\sigma x-y,\sigma^{2m}t)u_0(y)dy\\
&=\sigma^{\frac{2m-\alpha}{p-1}}\sigma^{n}\int_{\RR^{n}}g_{m}(\sigma x-\sigma y,\sigma^{2m}t)u_0(\sigma y)dy\\
&=\int_{\RR^{n}}g_{m}(x-y,t)\sigma^{\frac{2m-\alpha}{p-1}}u_0(\sigma y)dy\\
(E_{m,0}(t)u_0)_{\sigma}(x)&=E_{m,0}(t)u_{\sigma,0}(x).
\end{align*}
In a similar fashion we show that $(Tu(t))_{\sigma}=Tu_{\sigma}(t)$ so that  for any $\sigma>0$, $u_{\sigma}$ is a mild solution of Eq. \eqref{eq:main-eq}.
 On the other hand, if $u_0$ is homogeneous of degree $-\frac{2m-\alpha}{p-1}$, then we have  $\|u_0\|_{L^{q_c,\infty}}=\|u_{\sigma,0}\|_{L^{q_c,\infty}}$ and by hypothesis, $\|u_{\sigma,0}\|_{q_c,\infty}$ is small. An application of Theorem \ref{thm:GWP-Lp} yields the existence of $u_{\sigma}$ in $C_{b}([0,\infty);L^{q_c,\infty}(\RR^{n}))$	and by uniqueness (in a small closed ball) we deduce that $u(t)=u_{\sigma}(t)$ a.e. $x\in \RR^{n}$, $t>0$. The second assertion in Theorem \ref{thm:self-similarity} can be established in a similar way, the details are therefore omitted. 
\end{proof}

\begin{proof}[Proof of Corollary \ref{cor:trivial-self-similar}]
Assume that $u\in C_{b}([0,\infty);L^{q_c}(\RR^{n}))$ is a self-similar solution of Eq. \eqref{eq:mild-solution}, then $u(x,t)=\sigma^{\frac{2m-\alpha}{p-1}}u(\sigma x,\sigma^{2m}t)$, $t>0$ which by the change of variable $\sigma^{2m}t=1$ reads $u(x,t)=t^{-\frac{2m-\alpha}{2m(p-1)}}u(t^{-\frac{1}{2m}}x,1)$. Furthermore, we have  
\begin{align*}\displaystyle\lim_{t\rightarrow \infty}\|u(t)\|_{L^{q_c}}&=\displaystyle\lim_{t\rightarrow \infty}\|t^{-\frac{2m-\alpha}{2m(p-1)}}u(t^{-1/2m}\cdot,1)\|_{L^{q_c}}\\
&=\displaystyle\lim_{t\rightarrow \infty}t^{-\frac{2m-\alpha}{2m(p-1)}}t^{\frac{n}{2mq_c}}\|u(\cdot,1)\|_{L^{q_c}}\\
&=\|u(\cdot,1)\|_{L^{q_c}}.
\end{align*}
Since $u_0$ is small in $L^{q_c}(\RR^{n})$, Corollary \ref{cor:decay-of-sol} below tells us that $\displaystyle\lim_{t\rightarrow \infty}\|u(t)\|_{L^{q_c}}=0$. That is, $u(x,1)=0$ for a.e. $x\in \RR^{n}$ and hence $u=0$ a.e. on $\RR^{n}\times (0,\infty)$.
\end{proof}
The solution constructed in Theorem \ref{thm:GWP-Lp} enjoys other qualitative features which are inherent from the initial data. We collect some of those properties below.   
\begin{theorem}\label{thm:qualitative-features}
Let $u$ be the global-in-time solution given by Theorem \ref{thm:GWP-Lp}. One has:
\begin{enumerate}[label={(\roman*)}]
\item (Radial symmetry).\label{part:radialsymmetry} If $u_0$ is radial, then the solution $u$ is radial in $\RR^{n}$. 
\item (Radial monotonicity).\label{part:monotonicity} If $u_0$ is radially nonincreasing, then so is the solution $u$ in the spatial variable.
\item (Positivity).\label{part:positivity} Assume that $F$ preserves nonnegativity and let $m\in (0,1]$. If $u_0\geq 0$ and $u_0\neq 0$, then $u> 0$ a.e. $x\in \RR^{n}$, $t>0$. 	
\end{enumerate}
\end{theorem}
\begin{remark}Part \ref{part:positivity} of Theorem \ref{thm:qualitative-features} shows, in particular, that there are positive globally defined mild solutions of Eq. \eqref{eq:main-eq} whenever $m\in (0,1)$. This positivity, however, is exclusive to the latter case as far as one merely requires $u_0$ to be nonnegative. The failure of this property in the higher order case emanates from the fact that the kernel $g_m(\cdot,t)$ changes sign whenever $m$ is strictly larger than $1$. Thus, should positive global-in-time solutions exist, one may expect stronger assumptions on the data. In this direction, we refer the reader to the recent article \cite{Grunau} where the authors established existence of globally positive mild solutions of the semilinear and linear biharmonic heat equation (i.e. \eqref{eq:main-eq} in the case $m=2$ and $\alpha=0$) under certain assumptions on the initial data, prescribed in such a way that the changing sign phenomenon does not persist over large time. They further conjectured that their result is not restrictive to the biharmonic operator and should hold for the polyharmonic heat equation and its linear counterpart.  In parallel, it is of interest to investigate the eventual local positivity for solutions of Eq.  \eqref{eq:main-eq}, that is, positivity locally in time and space, see for instance \cite{Fer-Fer}. The techniques developed in the two previous references could be extended to also find necessary condition on the data giving rise to positive global-in-time solutions to Eq. \eqref{eq:main-eq}. We leave these questions for a future work.    
\end{remark}
\begin{proof}[Proof of Theorem \ref{thm:qualitative-features}]
As already mentioned, the solution $u$ obtained in Theorem \ref{thm:GWP-Lp} can be realized as the limit in a suitable sense of the sequence $(u_j)_{j\in \NN}$ defined as
\[u_1(\cdot,t)=E_{m,0}(t)u_0,\quad u_{j+1}(t)=\int_{0}^{t}E_{m,\alpha}(t-\tau)F(u_{j}(\tau))d\tau+u_1(t), \hspace{0.2cm} j\in \NN_{\geq 2}.
\] 	
Assume that $u_0$ is radial. As a convolution of $u_0$ with a radial kernel, $u_1(t)$ is radial in the spatial variable. By induction, we deduce that $u_j(t)$, $j=2,3,...$ is radial in $x\in \RR^{n}$ for all $t>0$ and converges in $L^{q_c}(\RR^{n})$ (if $u_0\in L^{q_c}(\RR^{n})$) or in $L^{q_c,\infty}(\RR^{n})$ (if $u_0\in L^{q_c,\infty}(\RR^{n})$). In either case, this convergence implies the existence of a subsequence of $(u_j)$ which converges almost everywhere to the solution $u$. Hence, \ref{part:radialsymmetry} immediately follows from the fact that almost everywhere convergence preserves radial symmetry. It is also known that the space of radial nonincreasing functions on $\RR^{n}$ is closed under convolution. Thus given that $g_m$ is radially nonincreasing, $u_1(t)$ for all $t>0$ is radially nonincreasing whenever $u_0$ has this property. Consequently, the sequence $(u_j(t))$ equally enjoys this feature by induction. The limit $u$, which is the solution of (\ref{eq:mild-solution}) is therefore radially nonincreasing. This proves \ref{part:monotonicity}. The last Part \ref{part:positivity} follows from the fact that $g_m$ has a positive kernel whenever $m\in (0,1]$.  
\end{proof}
\subsection{Long time asymptotic behavior}
We analyze the asymptotic stability of solutions in weak Marcinkiewicz space. It is shown that if two solutions of the linear equation $\partial_tu=-(-\Delta)^{m}u$ in $\RR^{n}\times (0,\infty)$ with data in $L^{q_c,\infty}(\RR^{n})$ are close for large times, then the corresponding solutions of the nonlinear equation remain close as $t\rightarrow \infty$. More precisely, we have:
\begin{theorem}\label{thm:stability}
Let $u_0,v_0$ be two small data in $L^{q_c,\infty}(\RR^{n})$ and consider $u,v$ two global solutions of Eq. \eqref{eq:mild-solution} with initial data $u_0$ and $v_0$, respectively. Assume that
$\displaystyle\lim_{t\rightarrow \infty}\|E_{m,0}(t)(u_0-v_0)\|_{L^{q_c,\infty}}=0.
$ Then 
\begin{equation}
\lim_{t\rightarrow\infty}\|u(t)-v(t)\|_{L^{q_c,\infty}}=0.
\end{equation}
Moreover, if $q>1$ and $\delta>0$ are as in \ref{GWP-P3}, then 
\begin{equation}
\lim_{t\rightarrow\infty}t^{\delta}\|u(t)-v(t)\|_{L^{q,\infty}}=0,
\end{equation}
whenever $\displaystyle\lim_{t\rightarrow \infty}t^{\delta}\|E_{m,0}(t)(u_0-v_0)\|_{L^{q,\infty}}=0.
$ 
\end{theorem}  
A close inspection shows that solutions in strong Lebesgue spaces have a rather easy-to-describe asymptotic behavior. To see this, let $u_0\in L^{q_c}(\RR^{n})$, by density there exists a sequence of functions $(u_{0N})_{N\in \NN}$ in $L^{q_c}\cap L^{q_0}(\RR^{n})$, $p_0\in (1,q_c)$ such that $u_{0N}\rightarrow u_0$ as $N\rightarrow \infty$ in $L^{q_c}(\RR^{n})$. Thus, for every $\varepsilon>0$ and $N$ large, we have  \[\|E_{m,0}(t)u_0\|_{L^{q_c}(\RR^{n})}\leq Ct^{-\frac{n}{2m}(\frac{1}{p_0}-\frac{1}{q_c})}\|u_{0N}\|_{L^{p_0}(\RR^n)}+c\varepsilon\]  
from which it follows that $\displaystyle\lim_{t\rightarrow \infty}\|E_{m,0}(t)u_0\|_{L^{q_c}}=0.$ By analogy,  $\displaystyle\lim_{t\rightarrow \infty}t^{\delta}\|E_{m,0}(t)u_0\|_{L^{q}}=0$. A direct consequence of these facts is that at large times, the solution $u$ decays to $0$ in both $L^{q_c}(\RR^{n})$ and the strong Lebesgue version of $S_q^{\delta}$, say $LS_q^{\delta}$. 
\begin{corollary}\label{cor:decay-of-sol}
Let $u\in C_{b}([0,\infty);L^{q_c}(\RR^{n}))$ (resp. $u\in LS_q^{\delta}$) be a global solution of Eq. \eqref{eq:mild-solution} subject to initial data $u_0\in L^{q_c}(\RR^{n})$ with small norm. Then 
\begin{equation}\label{eq:decay}
\lim_{t\rightarrow \infty}\|u(t)\|_{L^{q_c}}=0\quad(\text{resp}. \lim_{t\rightarrow \infty}t^{\delta}\|u(t)\|_{L^{q}}=0).
\end{equation} 
\end{corollary} 
 It is not hard to see that this asymptotic decay cannot hold in $L^{q_c,\infty}(\RR^{n})$ -- In fact, we have $\displaystyle\lim_{t\rightarrow \infty}\big\|E_{m,0}(t)|\cdot|^{-\frac{2m-\alpha}{p-1}}\big\|_{L^{q_c,\infty}}\neq 0.
 $  Yet, another consequence of Theorem \ref{thm:stability} is that a solution which initially is a suitable perturbation of a homogeneous function (of degree $-\frac{2m-\alpha}{(p-1)}$) in the critical weak Lebesgue space is asymptotic to a self-similar solution.
\begin{corollary}\label{cor:stability-perturbation}
Let $w_0\in L^{q_c,\infty}(\RR^{n})$ be sufficiently small and such that $w_0(\sigma \cdot)=\sigma^{-\frac{2m-\alpha}{(p-1)}}w_0$ for all $\sigma>0$ and call $w$ the global self-similar solution given by Theorem \ref{thm:self-similarity} such that $w(0)=w_0.$ Set $u_0=w_0+\varphi$ where $\varphi\in \overline{S(\RR^{n})}^{L^{q_c,\infty}(\RR^{n})}$ and assume $u_0$ is small. If $u$ is the global-in-time solution with initial data $u_0$, then 
\begin{equation}\label{eq:perturbation-eq}
\lim_{t\rightarrow \infty}\|u(t)-w(t)\|_{L^{q_c,\infty}}=0.
\end{equation}	
\end{corollary} 
Remark that $\displaystyle \lim_{t\rightarrow \infty}\|E_{m,0}(t)\varphi\|_{L^{q_c,\infty}(\RR^n)}=0$. To see this, observe that  by hypothesis, for each $\varepsilon>0$, there exists a sequence $(\varphi_n)$ in $S(\RR^n)$ with $\|\varphi_n-\varphi\|_{L^{q_c,\infty}(\RR^n)}<\varepsilon$.  By Proposition \ref{prop:smoothing-effect}, we find that
\begin{align*}
\|E_{m,0}(t)\varphi\|_{L^{q_c,\infty}(\RR^n)}\leq C\|E_{m,0}(t)\varphi_n\|_{L^{q_c,\infty}(\RR^n)}+C\varepsilon\leq Ct^{-\frac{n}{2m}(2-1/q_c)}\|\varphi_n\|_{L^2(\RR^n)}+C\varepsilon    
\end{align*}
since $S(\RR^n)\subset L^2(\RR^n)$ continuously. This being said, \eqref{eq:perturbation-eq} follows after applying Theorem \ref{thm:stability} which we now prove.  Observe in passing that from our previous discussion, choosing $v_0=0$ with a slight modification establishes Corollary \ref{cor:decay-of-sol}.
 \begin{proof}[Proof of Theorem \ref{thm:stability}]
 Let $u$ and $v$ two global solutions of (\ref{eq:mild-solution}) initially given by $u_0$ and $v_0$ respectively. Assume that $\displaystyle\lim_{t\rightarrow \infty}\|E_{m,0}(t)(u_0-v_0)\|_{L^{q_c,\infty}}=0$. For $t>0$, write 
 \[u(t)-v(t)=E_{m,0}(t)(u_0-v_0)+\int_{0}^{t}E_{m,\alpha}(t-s)[F(u(s))-F(v(s))]ds\] and put $L=\displaystyle\limsup_{t\rightarrow \infty}\|u(t)-v(t)\|_{L^{q_c,\infty}}$. Note that  $L\geq 0$ is finite. The goal is to show that $L=0$. One has 
 \begin{equation}\label{eq:mainest-stability}\|u(t)-v(t)\|_{L^{q_c,\infty}}\leq \|E_{m,0}(t)u_0-v_0\|_{L^{q_c,\infty}}+\bigg\|\int_{0}^{t}E_{m,\alpha}(t-s)[F(u(s))-F(v(s))]ds\bigg\|_{{L^{q_c,\infty}}}.
 \end{equation}
 Let $0<\mu<1$ (yet to be chosen) and set 
 \begin{equation*}
I^{\mu}=\bigg\|\int_{0}^{\mu t}E_{m,\alpha}(t-s)[F(u(s))-F(v(s))]ds\bigg\|_{_{L^{q_c,\infty}}}
 \end{equation*}
and 
\begin{equation*}
I_{\mu}=\bigg\|\int_{\mu t}^{t}E_{m,\alpha}(t-s)[F(u(s))-F(v(s))]ds\bigg\|_{{L^{q_c,\infty}}}.
\end{equation*}
With the aid of Proposition \ref{prop:smoothing-effect}, in addition to  $p>\frac{n-\alpha}{n-2m}$, we have 
\begin{align*}
I^{\mu}&\leq \int_{0}^{\mu t}\|E_{m,\alpha}(t-s)[F(u(s))-F(v(s))]\|_{{L^{q_c,\infty}}}ds\\
&\leq C\int_{0}^{\mu t}(t-s)^{-\frac{n(p-1)}{2mq_c}-\frac{\alpha}{2m}}\|[F(u(s))-F(v(s))]\|_{{L^{q_c,\infty}}}ds\\
&\leq C\int_{0}^{\mu t}(t-s)^{-\frac{n(p-1)}{2mq_c}-\frac{\alpha}{2m}}\|u(s)-v(s)\|_{{L^{q_c,\infty}}}(\|u(s)\|^{p-1}_{{L^{q_c,\infty}}}+\|v(s)\|^{p-1}_{{L^{q_c,\infty}}})ds\\
&\leq C2^{p}\vartheta^{p-1}\int_{0}^{\mu t}(t-s)^{-\frac{n(p-1)}{2mq_c}-\frac{\alpha}{2m}}\|u(s)-v(s)\|_{{L^{q_c,\infty}}}ds\\
&\leq C2^{p}\vartheta^{p-1}t^{-\frac{n(p-1)}{2mq_c}-\frac{\alpha}{2m}+1}\int_{0}^{\mu}(1-\tau)^{-\frac{n(p-1)}{2mq_c}-\frac{\alpha}{2m}}\|u(t\tau)-v(t\tau)\|_{{L^{q_c,\infty}}}d\tau\\
&\leq  C2^{p}\vartheta^{p-1}\int_{0}^{\mu}(1-\tau)^{-1}\|u(t\tau)-v(t\tau)\|_{{L^{q_c,\infty}}}d\tau
\end{align*}
where in the estimate before the last we have made the change of variable $s=\tau t$ and used the fact that $-\frac{n(p-1)}{2mq_c}-\frac{\alpha}{2m}+1=0$. Since $(1-\tau)^{-1}\in L^{1}_{loc}([0,\mu])$, it follows from the Dominated Convergence Theorem that 
\begin{equation}\label{eq:bound-on-Imu}
I^{\mu}\leq C2^{p}\vartheta^{p-1}\ln(1-\mu)^{-1}L.
\end{equation}
An obvious modification of the proof of Lemma \ref{lem:nonlinearbound} also yields the estimate 
\begin{equation}
I_{\mu}\leq C'2^{p}\vartheta^{p-1}\sup_{\tau\in [\mu t,t]}\|u(\tau)-v(\tau)\|_{{L^{q_c,\infty}}}
\end{equation}
so that in passing to the $\limsup$ as $t\rightarrow\infty$ on both sides of (\ref{eq:mainest-stability}) and taking (\ref{eq:bound-on-Imu}) into account, we arrive at 
\begin{equation*}
L\leq 2^{p}\vartheta^{p-1}(C\big|\ln(1-\mu)^{-1}\big|+C')L=\theta L.
\end{equation*}
Bearing in mind that $C^{\prime}2^{p}\vartheta^{p-1}<1$, one may choose $\mu\in (0,1)$ sufficiently small to allow $\theta<1$, and hence $L=0$. 
\end{proof}

\end{document}